\documentclass[aop,MSNbibl,seceqn,nameyear,dvips]{arximspdf}

\doi{10.1214/13-AOP899} 
\volume{43}
\issue{1}
\pubyear{2015}
\firstpage{240}
\lastpage{285}

\makeatletter
\newcommand{\rrvert}{\vert}
\newcommand{\llvert}{\vert}
\newcommand{\iint}{\mathop{\int\!\!\int}}
\newcommand{\bigintsss}{\int}
\newcommand{\idotsint}{\int\cdots\int}
\newtheorem{theorem}{Theorem}[section]
\newtheorem{lemma}[theorem]{Lemma}
\newtheorem{proposition}[theorem]{Proposition}
\newproclaim{remark}[theorem]{Remark}
\newproclaim{example}[theorem]{Example}
\newcommand{\SaS}{\mathrm{S}\alpha\mathrm{S}}
\makeatother

\begin{document}
\begin{frontmatter}

\title{Functional central limit theorem for heavy tailed stationary
infinitely divisible processes generated by conservative
flows\thanksref{T1}}
\runtitle{Functional central limit theorem}

\begin{aug}
\author[A]{\fnms{Takashi} \snm{Owada}\ead[label=e1]{to68@cornell.edu}}
\and
\author[B]{\fnms{Gennady} \snm{Samorodnitsky}\corref{}\ead[label=e2]{gs18@cornell.edu}}
\runauthor{T. Owada and G. Samorodnitsky}
\affiliation{Cornell University}
\address[A]{School of Operations Research\\
\quad and Information Engineering\\
Cornell University\\
Ithaca, New York 14853\\
USA\\
\printead{e1}} 
\address[B]{School of Operations Research\\
\quad and Information Engineering\\
Department of Statistical Science\\
Cornell University \\
Ithaca, New York 14853\\
USA\\
\printead{e2}}
\end{aug}
\thankstext{T1}{Supported in part by ARO Grants W911NF-07-1-0078
and W911NF-12-10385,
NSF Grant DMS-10-05903 and NSA Grant H98230-11-1-0154 at Cornell University.}

\received{\smonth{9} \syear{2012}}
\revised{\smonth{10} \syear{2013}}

%
\begin{abstract}
We establish a new class of functional central limit theorems for
partial sum of certain symmetric stationary infinitely divisible
processes with
regularly varying L\'evy measures. The limit process is a new class of
symmetric stable self-similar processes with stationary increments
that coincides on a part of its parameter space with a previously
described process. The normalizing sequence and the limiting process
are determined by the ergodic-theoretical properties of the flow
underlying the integral representation of the process. These
properties can be interpreted as determining how long the memory of
the stationary infinitely divisible process is. We also
establish functional convergence, in a strong distributional sense,
for conservative pointwise dual ergodic maps preserving an infinite
measure.
\end{abstract}

%

\end{frontmatter}

\section{Introduction} \label{secintro}
Let $\mathbf{X}=(X_1,X_2,\ldots)$ be a discrete time stationary\break stochastic
process. A (functional) central limit theorem for such a process is a
statement of the type
%
%
\begin{equation}
\label{eFCLT} \Biggl( \frac{1}{c_n}\sum_{k=1}^{\lceil nt\rceil}
X_k -h_nt, 0\leq t\leq1 \Biggr) \Rightarrow \bigl( Y(t),
0\leq t\leq1 \bigr).
\end{equation}
Here, $(c_n)$ is a positive sequence growing to infinity, $(h_n)$ a
real sequence, and $ ( Y(t),   0\leq t\leq1 )$ is a
nondegenerate (i.e., nondeterministic) process. Convergence in
(\ref{eFCLT}) is at least in finite-dimensional distributions, but
preferably it is a weak convergence in the space $D[0,1]$ equipped
with an appropriate topology. Not every stochastic process satisfies a
central limit theorem, and for those that do, it is well known that
both the rate of growth of the scaling constant $c_n$ and the nature
of the limiting process $\mathbf{Y}= ( Y(t),   0\leq t\leq
1 )$
are determined both by the marginal tails of the stationary process
$\mathbf{X}$ and its dependence structure. The limiting process (under very
minor assumptions) is necessarily self-similar with stationary
increments; this is known as the Lamperti theorem; see
\citet{lamperti1962}.

If, say, $X_1$ has a finite second
moment, and $\mathbf{X}$ is an i.i.d. sequence then, clearly, one can choose
$c_n=n^{1/2}$, and then $\mathbf{Y}$ is a Brownian motion. With
equally light
marginal tails, if the memory is sufficiently short, then one expects the
situation to remain, basically, the same, and this turns out to be the
case. When the variance is finite, the basic tool to measure
dependence is, obviously, the correlations, which have to decay fast
enough. It is well known, however, that a fast decay of correlations
is alone not sufficient for this purpose, and, in general, certain
strong mixing conditions have to be assumed. See, for example,
\citet{rosenblatt1956} and, more recently,
\citet{merlevedepeligradutev2006}. If the memory is not sufficiently
short, then both the rate of growth of $c_n$ can be different from
$n^{1/2}$, and the limiting process can be different from the Brownian
motion. In fact, the limiting process may fail to be Gaussian at all;
see, for example, \citet{dobrushinmajor1979} and \citet{taqqu1979}.

If the marginal tails of the process are heavy, which in this case
means that $X_1$ is in the domain of attraction of an $\alpha$-stable
law, $0<\alpha<2$, and $\mathbf{X}$ is an i.i.d. sequence then
clearly one
can choose $c_n$ to be the inverse of the marginal tail (this makes
$c_n$ vary regularly with exponent $1/\alpha$), and then $\mathbf{Y}$
is an
$\alpha$-stable L\'evy motion. Again, one expects the situation to
remain similar if the memory is sufficiently short. Since correlations
do not exist under heavy tails, statements of this type have been
established for special models, often for moving average models; see,
for example, \citet{davisresnick1985}, \citet{avramtaqqu1992} and
\citet{paulauskassurgailis2008}. Once again, as the memory gets
longer, then both the rate of growth of $c_n$ can be different from
that obtained by inverting the marginal tail, and the limiting process
will no longer have independent increments (i.e., be an
$\alpha$-stable L\'evy motion). It is here, however, that the picture
gets more interesting than in the case of light tails. First of all,
in absence of correlations there is no canonical way of measuring how
much longer the memory gets. Even more importantly, certain types of
memory turn out to result in the limiting process $\mathbf{Y}$ being a
self-similar $\alpha$-stable process with stationary increments of a
canonical form, the so-called linear fractional stable motion; see,
for example, \citet{maejima1983B} for an example of such a situation, and
\citet{samorodnitskytaqqu1994} for information on self-similar
processes. However, when the memory gets even longer, linear
fractional stable motions disappear as well, and even more ``unusual''
limiting processes $\mathbf{Y}$ may appear. This phenomenon may
qualify as
change from short to long memory; see \citet{samorodnitsky2006LRD}.

In this paper, we consider a functional central limit theorem for a
class of heavy tailed stationary processes exhibiting long memory in
this sense. It is particularly interesting both because of the manner
in which
memory in the process is measured, and because the limiting process
$\mathbf{Y}$ that happens to be an extension of a very recently discovered
self-similar stable process with stationary increments. Specifically, we
will assume that $\mathbf{X}$ is a stationary infinitely divisible process (satisfying
certain assumptions, described in detail in Section~\ref{secprelim}).
That is, all finite-dimensional distributions of
$\mathbf{X}$ are infinitely divisible; we refer the reader to \citet{rajputrosinski1989} for
more information on infinitely divisible processes and their integral
representations
we will work with in the sequel.

The class of central limit theorems we consider involves a significant
interaction of probabilistic and ergodic-theoretical ideas and
tools. To make the discussion more transparent, we will only
consider symmetric infinitely divisible processes without a Gaussian
component (but
there is no doubt that results of this type will hold in a greater
generality as well). The law of such a process is
determined by its (function level) L\'evy measure. This is a (uniquely
determined) symmetric measure $\kappa$ on $\mathbb{R}^\mathbb{N}$ satisfying
\[
\kappa \bigl( {\mathbf x}=(x_1,x_2,\ldots)\in
\mathbb{R}^\mathbb{N}\dvtx  x_j=0\mbox{ for all }j\in\mathbb{N}
\bigr)=0
\]
and
\[
\int_{\mathbb{R}^\mathbb{N}}\min\bigl(1,x_j^2\bigr)
\kappa(d{\mathbf x})<\infty\qquad \mbox{for each }j\in\mathbb{N},
\]
such that for each finite subset $\{ j_1,\ldots, j_k\}$ of $\mathbb
{N}$, the
$k$-dimensional L\'evy measure of the infinitely divisible random
vector $ (
X_{j_1}, \ldots, X_{j_k} )$ is given by the projection of $\kappa
$ on
the appropriate coordinates of ${\mathbf x}$; see \citet{maruyama1970}.

Because of the stationarity of the process $\mathbf{X}$, its L\'evy measure
$\mu$ is invariant under the left shift $\theta$ on $\mathbb
{R}^\mathbb{N}$,
\[
\theta(x_1,x_2,x_3,\ldots) =
(x_2,x_3,\ldots).
\]
It has been noticed in the last
several years that the ergodic-theoretical properties of the shift
operator with respect to the L\'evy measure have a profound effect on
the memory of the stationary process $\mathbf{X}$. The L\'evy measure
of the
process is often described via an integral representation of the
process, and in some cases the shift operator with respect to the
L\'evy measure can be related to an operator acting on the space on
which the integrals are taken. Thus,
\citet{rosinskisamorodnitsky1996} and \citet{samorodnitsky2005a}
dealt with the ergodicity and mixing of stationary stable processes,
while \citet{roy2008} dealt with general stationary infinitely
divisible processes. The effect of the ergodic-theoretical properties of the shift
operator with respect to the L\'evy measure on the partial maxima of
stationary stable processes was discussed in
\citet{samorodnitsky2004a}.

In the present paper, we consider stationary symmetric infinitely
divisible processes
without a Gaussian component given via an integral representation
described in Section~\ref{secprelim}. This representation naturally
includes a measure-preserving operator on a measurable space, and we
related its ergodic-theoretical properties to the kind of central
limit theorem the process satisfies. We consider the so-called
conservative operators that turn out to lead to nonstandard limit
theorems of the type that, to the best of our knowledge, have not been
observed before.

We describe our setup in Section~\ref{secprelim}. In Section~\ref{seclimprocess}, we introduce the limiting symmetric
$\alpha$-stable (henceforth, $\SaS$) self-similar
process with stationary increments and discuss its
properties. In Section~\ref{secergodic}, we present the
ergodic-theoretical notions that we use in the paper.
The exact assumptions in the central limit theorem are stated in
Section~\ref{secCLT}. In this section, we also present the statement
of the theorem and several examples. The proof of the theorem uses
several distributional ergodic-theoretical results we present and
prove in Section~\ref{secdistrresults}. These results may be of
independent interest in ergodic theory. Finally, the proof of the
central limit theorem is completed in Section~\ref{secproof}.

\section{The setup} \label{secprelim}
We consider infinitely divisible processes of the form
%
%
\begin{equation}
\label{etheprocess} X_n = \int_E
f_n(x) \,dM(x), \qquad n=1,2,\ldots,
\end{equation}
where $M$ is an infinitely divisible random measure on a measurable space
$(E,\mathcal{E})$, and the functions $f_n,  n=1,2,\ldots$ are
deterministic functions of the form
%
%
\begin{equation}
\label{ethekernel} f_n(x) = f\circ T^n(x) = f \bigl(
T^nx \bigr),\qquad x\in E, n=1,2,\ldots,
\end{equation}
where $f\dvtx E \to\mathbb{R}$ is a measurable function, and $T\dvtx E \to E$
a measurable map. The (independently scattered) infinitely divisible random measure
$M$ is assumed to be a homogeneous symmetric infinitely divisible random measure
without a Gaussian component, with control measure $\mu$ and local
L\'evy measure $\rho$. That is, $\mu$ is a \mbox{$\sigma$-}finite measure on
$E$, which we will assume to be infinite. Further, $\rho$ is a
symmetric L\'evy measure on $\mathbb{R}$, and for every $A \in
\mathcal{E}$
with $\mu(A) < \infty$, $M(A)$ is a (symmetric) infinitely
divisible random variable such that
%
%
\begin{equation}
\label{eqchfIDRM} E e^{iu M(A)} = \exp \biggl\{ -\mu(A) \int
_{\mathbb{R}} \bigl(1-\cos (ux) \bigr) \rho(dx) \biggr\}, \qquad u \in
\mathbb{R}.
\end{equation}

It is clear that, in order for the process $\mathbf{X}$ to be well defined,
the functions $f_n,  n=1,2,\ldots$ have to satisfy certain
integrability assumptions; the assumptions we will impose below will
be sufficient for that. Once the process $\mathbf{X}$ is well defined,
it is,
automatically, symmetric and infinitely divisible, without a Gaussian
component, with
the function level L\'evy measure given by
%
%
\begin{equation}
\label{emeasureprocess} \kappa= (\rho\times\mu) \circ K^{-1}
\end{equation}
with $K\dvtx   \mathbb{R}\times E \to\mathbb{R}^\mathbb{N}$ given by
$K(x,s) = x (
f_1(s),f_2(s),\ldots )$, $s\in E,  x\in\mathbb{R}$. For
details, see
\citet{rajputrosinski1989}.

We will assume that the measurable map $T$ preserves the control
measure $\mu$. It follows immediately from (\ref{emeasureprocess})
and the form of the functions $(f_n)$ given in (\ref{ethekernel})
that the L\'evy measure $\kappa$ is invariant under the left shift
$\theta$, and hence, the process $\mathbf{X}$ is stationary. We
intend to
relate the ergodic-theoretical properties of the map $T$ to the
dependence properties of the process $\mathbf{X}$, and subsequently,
to the
kind of central limit theorem the process satisfies. We refer the
reader to \citet{aaronson1997} for more details on the
ergodic-theoretical notions used in the sequel. A~short review of what
we need will be given in Section~\ref{secergodic} below.

Our basic assumption is that the map $T$ is conservative. This
property has already been observed to be related to long memory in the
process $\mathbf{X}$; see, for example, \citet{samorodnitsky2004a} and
\citet{roy2008}. We will quantify the resulting length of memory by
assuming further that the map $T$ is ergodic and pointwise dual
ergodic, with a regularly varying normalizing sequence. We will see
that the
exponent of regular variation plays a major role in the central limit
theorem.

The second major ``player'' in the central limit theorem is
the heaviness of the marginal tail of the process $\mathbf{X}$. We
will assume
that the local L\'evy measure $\rho$ has a regularly varying tail with
index $-\alpha$, $0 < \alpha< 2$, that is,
%
%
\begin{equation}
\label{eregvarlevy} \rho(\cdot,\infty) \in RV_{-\alpha}\qquad\mbox{at infinity.}
\end{equation}
With a proper integrability assumption on the function $f$ in
(\ref{ethekernel}), the process $\mathbf{X}$ has regularly varying marginal
(and even finite-dimensional) distributions, with the same tail exponent
$-\alpha$; see \citet{rosinskisamorodnitsky1993}. That is, all the
finite-dimensional distributions of the process are in the domain of
attraction of a $\SaS$ law.

This leads to a rather satisfying picture, in which the kind of the
central limit theorem that holds for the process $\mathbf{X}$ depends
both on the marginal tails of the process and on the length of memory
in it, and both are clearly parameterized.

In fact, in order to obtain the central limit theorem for the process
$\mathbf{X}$, we will need to impose more specific assumptions on the map
$T$. We will also, clearly, need specific integrability assumptions on
the kernel in the integral representation of the process. These
assumptions are presented in Section~\ref{secCLT}.

We proceed, first, with a description of the limiting process we will
eventually obtain.

\section{The limiting process} \label{seclimprocess}

In this section, we will introduce a class of self-similar $\SaS$ processes with stationary increments. These processes will later
appear as weak limits in the central limit theorem. We will see this
process is an extension (to a wider range of parameters) of a class
recently introduced by \citet{dombryguillotin-plantard2009}. Before
introducing this process, we need to do some preliminary work.

For $0<\beta<1$, let $ ( S_{\beta}(t), t\geq0 )$ be a
\mbox{$\beta$-}stable subordinator, that is, a L\'evy process with increasing
sample paths, satisfying $Ee^{-\theta S_{\beta}(t)} = \exp\{ -t
\theta^{\beta} \}$ for \mbox{$\theta\geq0$ and $t\geq0$}; see,
for example, Chapter III of \citet{bertoin1996}. Define its inverse
process by
%
%
\begin{equation}
\label{eMLprocess} M_{\beta}(t) = S_{\beta}^{\leftarrow}(t) = \inf
\bigl\{u\geq0\dvtx  S_{\beta}(u) \geq t \bigr\},\qquad t\geq0.
\end{equation}
Recall that the marginal distributions of the process $ (
M_{\beta}(t),   t\geq0 )$ are the Mittag--Leffler distributions,
with the Laplace transform
%
%
\begin{equation}
\label{MTtransform} E \exp\bigl\{ \theta M_{\beta}(t) \bigr\} = \sum
_{n=0}^{\infty} \frac{
(\theta t^{\beta})^n}{\Gamma(1+n \beta)}, \qquad\theta\in
\mathbb{R};
\end{equation}
see Proposition 1(a) in \citet{bingham1971}. We will call this process
\textit{the Mittag--Leffler process}. This process has a continuous and
nondecreasing version; we will always assume that we are working with
such a version. It follows from (\ref{MTtransform}) (or simply from
the definition) that the Mittag--Leffler process is self-similar with
exponent $\beta$. Further, all of its moments are finite. Recall,
however, that this process has neither stationary nor independent
increments; see, for example, \citet{meerschaertscheffler2004}.

We are now ready to introduce the new class of self-similar $\SaS$ processes with stationary increments announced at the beginning of
this section. Let $0<\alpha<2$ and $0<\beta<1$, and let
$(\Omega^{\prime},\mathcal{F}^{\prime},P^{\prime})$ be a probability
space. We define
%
%
\begin{equation}
\label{eqMLSM} Y_{\alpha,\beta}(t) = \int_{\Omega^{\prime} \times[0,\infty)}
M_{\beta} \bigl((t-x)_+,\omega^{\prime} \bigr) \,d Z_{\alpha,\beta}
\bigl(\omega^{\prime},x\bigr), \qquad t \geq0,
\end{equation}
where $Z_{\alpha,\beta}$ is a $\SaS$ random measure on $\Omega
^{\prime}
\times[0,\infty)$ with control measure $P^{\prime} \times\nu$, with
$\nu$ a measure on $[0,\infty)$ given by
$\nu(dx) = (1-\beta) x^{-\beta}   \,dx,   x>0$. Here, $M_\beta$~is a
Mittag--Leffler process defined on
$(\Omega^{\prime},\mathcal{F}^{\prime},P^{\prime})$. The random measure
$Z_{\alpha,\beta}$ itself, and hence, also the process
$Y_{\alpha,\beta}$, are defined on some generic probability space
$(\Omega,\mathcal{F},P )$. We refer the reader to
\citet{samorodnitskytaqqu1994} for more information on integrals with
respect to stable random measures.

In Theorem~\ref{tMTFSNbasic} below, we prove that the process
$ ( Y_{\alpha,\beta}(t),  t\geq0 )$ is a well-defined
self-similar $\SaS$ processes with stationary increments. We call it \textit{the
\mbox{\mbox{$\beta$-}}Mittag--Leffler} (or \mbox{$\beta$-}ML) \textit{fractional $\SaS$ motion}.

%
\begin{theorem} \label{tMTFSNbasic}
The \mbox{$\beta$-}ML fractional $\SaS$  motion is a well-defined self-similar
$\SaS$ processes with stationary increments. It is also self-similar with
exponent of self-similarity $H=\beta+ (1-\beta)/\alpha$.
\end{theorem}

\begin{pf}
By the monotonicity of the process $M_{\beta}$ we have, for any $t\geq
0$,
\[
\int_{[0,\infty)} \int_{\Omega^{\prime}} M_{\beta}
\bigl((t-x)_+,\omega^{\prime}\bigr)^{\alpha} P^{\prime}(d\omega)
\nu(dx)\leq t^{1-\beta} E^{\prime} M_{\beta}(t)^{\alpha}
< \infty,
\]
which proves that the process $ ( Y_{\alpha,\beta}(t),  t\geq
0 )$ is well defined. Further, by the \mbox{$\beta$-}self-similarity of
the process $M_{\beta}$, we have for any $k\geq1$, $t_1, \ldots, t_k
\geq0$, and $c>0$, for all real $\theta_1, \ldots, \theta_k $,
\begin{eqnarray*}
&& E \exp \Biggl\{ i \sum_{j=1}^k
\theta_j Y_{\alpha,\beta}(ct_j) \Biggr\}
\\
&&\qquad = \exp \Biggl
\{ - \int_0^{\infty} E^{\prime} \Biggl|\sum
_{j=1}^k \theta_j M_{\beta}
\bigl((ct_j-x)_+\bigr) \Biggr|^{\alpha} (1-\beta)x^{-\beta} \,dx
\Biggr\}
\\
&&\qquad = \exp \Biggl\{ - \int_0^{\infty} E^{\prime}
\Biggl|\sum_{j=1}^k \theta_j
c^H M_{\beta}\bigl((t_j-y)_+\bigr)
\Biggr|^{\alpha} (1-\beta)y^{-\beta} \,dy \Biggr\}
\\
&&\qquad = E \exp \Biggl\{ i \sum
_{j=1}^k \theta_j
c^H Y_{\alpha,\beta}(t_j) \Biggr\},
\end{eqnarray*}
which shows the $H$-self-similarity of the \mbox{$\beta$-}ML fractional $\SaS$
motion.

For the proof of stationary increment property, it suffices to check that
\[
E \exp \Biggl\{ i \sum_{j=1}^k
\theta_j \bigl(Y_{\alpha,\beta
}(t_j+s) -
Y_{\alpha,\beta}(s) \bigr) \Biggr\} = E \exp \Biggl\{ i \sum
_{j=1}^k \theta_j Y_{\alpha,\beta}(t_j)
\Biggr\}
\]
for all $k\geq1$, $t_1, \ldots, t_k \geq0$, $s\geq0$, and $\theta_1,
\ldots,\theta_k \in\mathbb{R}$. This is equivalent to verifying the
equality in
\begin{eqnarray*}
&& \int_0^{\infty} E^{\prime} \Biggl|\sum
_{j=1}^k \theta_j \bigl\{
M_{\beta}\bigl((t_j+s-x)_+\bigr) - M_{\beta}
\bigl((s-x)_+\bigr)\bigr\} \Biggr|^{\alpha
}x^{-\beta} \,dx
\\
&&\qquad =\int_0^{\infty} E^{\prime} \Biggl|\sum
_{j=1}^k \theta_j M_{\beta}
\bigl((t_j-x)_+\bigr) \Biggr|^{\alpha} x^{-\beta} \,dx.
\end{eqnarray*}
Changing variable by $r=s-x$ in the left-hand side and rearranging the
terms shows that we need to check the equality in
%
%
\begin{eqnarray}\label{esicheck1}
&& \int_{0}^s E^{\prime} \Biggl|
\sum_{j=1}^k \theta_j
\bigl(M_{\beta
}(t_j+r) - M_{\beta}(r)\bigr)
\Biggr|^{\alpha} (s-r)^{-\beta} \,dr
\nonumber\\[-8pt]\\[-8pt]
&&\qquad = \int_0^{\infty} E^{\prime} \Biggl|\sum
_{j=1}^k \theta_j M_{\beta}
\bigl((t_j-x)_+\bigr) \Biggr|^{\alpha} \bigl(x^{-\beta} -
(s+x)^{-\beta}\bigr) \,dx.\nonumber
\end{eqnarray}

Let $\delta_r = S_\beta ( M_\beta(r) )-r$ be the overshoot of
the level $r>0$ by the \mbox{$\beta$-}stable subordinator $ (
S_{\beta}(t), t\geq0 )$ related to $ (
M_{\beta}(t), t\geq0 )$ by (\ref{eMLprocess}). The law of
$\delta_r$ is known to be given by
%
%
\begin{equation}
\label{eqDynkinLamperti} P(\delta_r \in dx) = \frac{\sin\beta\pi}{\pi}
r^{\beta} (r+x)^{-1} x^{-\beta} \,dx,\qquad x>0;
\end{equation}
see, for example, Exercise 5.6 in \citet{kyprianou2006}. Further, by
the strong
Markov property of the stable subordinator we have
\[
\bigl( S_\beta \bigl( M_{\beta}(r)+t \bigr), t\geq0 \bigr)
\stackrel{d} {=} \bigl( r+\delta_r + S_\beta(t), t\geq0
\bigr),
\]
where $S_\beta$ and $\delta_r$ in the right-hand
side are independent. Therefore,
\begin{eqnarray*}
&& \bigl( M_{\beta}(t+r) - M_{\beta}(r), t\geq0 \bigr)
\\
&&\qquad = \bigl( \inf
\bigl\{u\geq0\dvtx  S_{\beta} \bigl( M_{\beta
}(r)+u \bigr)\geq t +r
\bigr\}, t\geq0 \bigr)
\\
&&\qquad \stackrel{d} {=} \bigl( \inf \bigl\{u\geq0\dvtx  S_{\beta}(u)\geq t-\delta
_r \bigr\}, t\geq0 \bigr)
\\
&&\qquad = \bigl( M_{\beta}\bigl((t-
\delta_r)_+\bigr), t\geq0 \bigr);
\end{eqnarray*}
once again, $M_\beta$ and $\delta_r$ in the right-hand
side are independent. We conclude that
%
%
\begin{eqnarray}\label{esicheck2}
&& \int_{0}^s E^{\prime} \Biggl|
\sum_{j=1}^k \theta_j
\bigl(M_{\beta}(t_j+r) - M_{\beta}(r)
\bigr)\Biggr|^{\alpha} (s-r)^{-\beta} \,dr\nonumber
\\
&&\qquad = \frac{\sin\beta\pi}{\pi} \int_0^{\infty}\!\! \int
_0^s E^{\prime} \Biggl|\sum
_{j=1}^k \theta_j M_{\beta}
\bigl((t_j-x)_+\bigr)\Biggr|^{\alpha}
\\
&&\hspace*{96pt}{}\times  r^{\beta}
(r+x)^{-1} x^{-\beta} (s-r)^{-\beta} \,dr \,dx.\nonumber
\end{eqnarray}
Using the integration formula,
\[
\int_0^1 \biggl(\frac{t}{1-t}
\biggr)^\beta\frac{1}{t+y} \,dt = \frac{\pi}{\sin\beta\pi} \biggl[ 1- \biggl(
\frac{y}{1+y} \biggr)^\beta \biggr],\qquad y>0,
\]
given on page 338 of \citet{gradshteynryzhik1994}, shows that
(\ref{esicheck2}) is equivalent to (\ref{esicheck1}). This
completes the proof.
\end{pf}

Recall that, when $0<\beta\leq1/2$, the Mittag--Leffler process of
(\ref{eMLprocess}) is distributionally equivalent to the local time
at zero of a symmetric stable L\'evy process with index of stability
$\hat\beta=(1-\beta)^{-1}$. Specifically,\vspace*{-1pt} let $ (
W_{\hat\beta}(t),   t\geq0 )$ be a symmetric $\hat\beta$-stable
L\'evy process, such that $Ee^{ir W_{\hat\beta}(t)} = \exp\{ -t
|r|^{\hat\beta} \}$ for $r\in\mathbb{R}$ and $t\geq0$. This
process has a
jointly continuous local time process, $L_t(x),   t\geq0,
x\in\mathbb{R}$; see, for example, \citet{getoorkesten1972}. Then
%
%
\begin{equation}
\label{eMLLT} \bigl( M_{\beta}(t), t\geq0 \bigr) \stackrel{d} {=} \bigl(
c_\beta L_t(0), t\geq0 \bigr)
\end{equation}
for some $c_{\beta} > 0$; see Section~11.1.1 in
\citet{marcusrosen2006}. Therefore, in the range $0<\beta\leq1/2$,
the \mbox{$\beta$-}ML fractional $\SaS$  motion (\ref{eqMLSM}) can be
represented in law as
%
%
\begin{equation}
\label{eqMLSM1} Y_{\alpha,\beta}(t) = c_\beta\int_{\Omega^{\prime} \times
[0,\infty)}
L_{(t-x)_+} \bigl( 0,\omega^{\prime} \bigr) \,d Z_{\alpha,\beta}\bigl(
\omega^{\prime},x\bigr), \qquad t \geq0,
\end{equation}
where $ ( L_t(x) )$ is the local time of a symmetric $\hat
\beta$-stable
L\'evy process defined on $(\Omega^{\prime},\mathcal{F}^{\prime
},P^{\prime})$.
Recall also
the $\hat\beta$-stable local time fractional $\SaS$  motion introduced
in \citet{dombryguillotin-plantard2009} [see also
\citet{cohensamorodnitsky2006}]. That process can be defined by
%
%
\begin{equation}
\label{eqLTSM} \widehat Y_{\alpha,\beta}(t) = \int_{\Omega^{\prime} \times\mathbb{R}}
L_{t} \bigl( x,\omega^{\prime} \bigr) \,d \widehat Z_{\alpha}
\bigl(\omega^{\prime},x\bigr), \qquad t \geq0,
\end{equation}
where $\widehat Z_{\alpha}$ is a $\SaS$ random measure on
$\Omega^{\prime}
\times\mathbb{R}$ with control measure $P^{\prime} \times\mathrm{Leb}$. We
claim that, in fact, if $0<\beta\leq1/2$,
%
%
\begin{equation}
\label{eequalSM} \bigl( Y_{\alpha,\beta}(t) t\geq0 \bigr) \stackrel {d}
{=}c_\beta ^{(1)} \bigl( \widehat Y_{\alpha,\beta}(t) t\geq0
\bigr)
\end{equation}
for some multiplicative constant $c_\beta^{(1)}$. Therefore,\vspace*{-1pt} one can
view the ML fractional~$\SaS$  motion as an extension of the
$\hat\beta$-stable local time fractional $\SaS$  motion from the range
$1<\hat\beta\leq2$ to the range $1<\hat\beta<\infty$. It is
interesting to note that the central limit theorem in Section~\ref
{secCLT} is of a very different type from the random walk in random
scenery situation of \citet{cohensamorodnitsky2006} and
\citet{dombryguillotin-plantard2009}.

To check (\ref{eequalSM}), let
\[
H_x=\inf \bigl\{ t\geq0\dvtx  W_{\hat\beta}(t)=x \bigr\},\qquad x\in
\mathbb{R}.
\]
Since $1<\hat\beta\leq2$, $H_x$ is a.s. finite for any $x\in\mathbb{R}$;
see, for example, Remark 43.12 in \citet{sato1999}. Further, by the strong
Markov property, for every $x\in\mathbb{R}$, the conditional law of $
( L_{H_x+t}(x),   t\geq0 )$ given
$\mathcal{F}^{\prime}_{H_x}$, coincides a.s. with the law of $  (
L_{t}(0),
t\geq0 ) $. We conclude that for any $k\geq1$, $t_1, \ldots, t_k
\geq0$, and real $\theta_1, \ldots, \theta_k $,
\begin{eqnarray*}
-\log E\exp \Biggl\{ \sum_{j=1}^k
\theta_j \widehat Y_{\alpha,\beta}(t_j) \Biggr\} &=& \int
_\mathbb{R}E^\prime \Biggl| \sum
_{j=1}^k \theta_j L_{t_j}(x)
\Biggr|^\alpha \,dx
\\
&=& \int_\mathbb{R}\int_0^\infty
E^\prime \Biggl| \sum_{j=1}^k
\theta_j L_{(t_j-y)_+}(0) \Biggr|^\alpha F_x(dy)
\,dx,
\end{eqnarray*}
where $F_x$ is the law of $H_x$. Using the obvious fact that
$H_x\stackrel{d}{=}
|x|^{\hat\beta}H_1$, an easy calculation shows that the mixture
$\int_\mathbb{R}F_x  \,dx$ is, up to a multiplicative constant, equal
to the
measure $\nu$ in (\ref{eqMLSM}). Therefore, for some constant
$c_\beta^{(1)}$,
\[
-\log E\exp \Biggl\{ \sum_{j=1}^k
\theta_j c_\beta^{(1)}\widehat Y_{\alpha,\beta}(t_j)
\Biggr\} = -\log E\exp \Biggl\{ \sum_{j=1}^k
\theta_j Y_{\alpha,\beta}(t_j) \Biggr\}
\]
and (\ref{eequalSM}) follows.

%
\begin{remark}\label{rkrangeH}
It is interesting to observe that, for a fixed $0<\alpha<2$, the range
of the exponent of self-similarity $H=\beta+ (1-\beta)/\alpha$ of the
\mbox{$\beta$-}ML fractional $\SaS$  motion, as
$\beta$ varies between 0 and 1, is a proper subset of the feasible
range of the exponent of self-similarity of stationary increment
self-similar $\SaS$ processes, which is $0 < H \leq\max(1, 1/\alpha)
$; see \citet{samorodnitskytaqqu1994}.
\end{remark}

It was shown in \citet{dombryguillotin-plantard2009} that the stable
local time fractional $\SaS$  motion is H\"older continuous. We extend
this statement to the ML fractional $\SaS$  motion.

%
\begin{theorem} \label{tholder}
The \mbox{$\beta$-}ML fractional $\SaS$  motion satisfies, with probability 1,
\[
\sup_{0\leq s<t\leq1/2}\frac{ |
Y_{\alpha,\beta}(t)-Y_{\alpha,\beta}(s) |}{(t-s)^\beta
|\log(t-s) |^{1-\beta}} <\infty
\]
if $0<\alpha<1$, and
\[
\sup_{0\leq s<t\leq1/2}\frac{ |
Y_{\alpha,\beta}(t)-Y_{\alpha,\beta}(s) |}{(t-s)^\beta
|\log(t-s) |^{3/2-\beta}} <\infty
\]
if $1\leq\alpha<2$.
\end{theorem}

\begin{pf}
The statement of the theorem follows from Lemma~\ref{lMLholder} and
the argument in Theorem 5.1 in \citet{cohensamorodnitsky2006}; see
also Theorem 1.5 in \citet{dombryguillotin-plantard2009}.
\end{pf}

The next lemma establishes H\"older continuity of the Mittag--Leffler
process~(\ref{eMLprocess}). The statement might be known, but we
could not find a reference, so we present a simple argument. In the
case $0<\beta\leq1/2$ (most of) the statement is in Theorem 2.1 in
\citet{ehm1981}, through the relation with the local time
(\ref{eMLLT}).

%
\begin{lemma} \label{lMLholder}
For $B>0$, let
\[
K = \sup_{0\leq s<t<s+1/2\leq B}\frac{ |
M_{\beta}(t)-M_{\beta}(s) |}{(t-s)^\beta
|\log(t-s) |^{1-\beta}}.
\]
Then $K$ is an a.s. finite random variable with all finite moments.
\end{lemma}

\begin{pf}
Because of the self-similarity of the Mittag--Leffler
process, it is enough to consider $B=1/2$. In the
course of the proof, we will use the notation $c(\beta)$ for a finite
positive constant that may depend on $\beta$, and that may change
from one appearance to another. Recall the lower tail estimate of a positive
\mbox{$\beta$-}stable random variable:
%
%
\begin{equation}
\label{elowertail} P \bigl( S_\beta(1)\leq\theta \bigr) \leq\exp \bigl\{
-c(\beta) \theta^{-\beta/(1-\beta)} \bigr\},\qquad 0<\theta\leq1;
\end{equation}
see \citet{zolotarev1986}. Let $\lambda\geq1$. We have
\begin{eqnarray*}
P(K>\lambda) &\leq& \sum_{n=1}^\infty P \Bigl(
\mathop{\sup_{0\leq s<t\leq1/2}}_{2^{-(n+1)}\leq t-s\leq2^{-n}} M_{\beta}(t)-M_{\beta}(s)>c(\beta)
\lambda n^{1-\beta}2^{-n\beta} \Bigr)
\\
&:=& \sum _{n=1}^\infty q_n(\lambda).
\end{eqnarray*}
For $n=1,2,\ldots,$ we use the following decomposition:
\begin{eqnarray*}
q_n(\lambda)&\leq& P \bigl( S_\beta(\lambda\log n)\leq1/2
\bigr)
\\
&&{} + P \bigl[\mbox{for some } 0<t\leq\lambda\log n, S_\beta \bigl( t+c(
\beta)\lambda n^{1-\beta}2^{-n\beta} \bigr) - S_\beta(t)
\leq2^{-n} \bigr]
\\
&:=& q ^{(1)}_n(
\lambda)+q^{(2)}_n(\lambda).
\end{eqnarray*}
Using (\ref{elowertail}) and self-similarity of the stable
subordinator, we obtain
\[
\sum_{n=1}^\infty q ^{(1)}_n(
\lambda) \leq c(\beta)^{-1}\exp \bigl\{ -c(\beta)\lambda^{1/(1-\beta)}
\bigr\}.
\]
On the other hand,
\begin{eqnarray*}
q^{(2)}_n(\lambda) &\leq& P \bigl( S_\beta
\bigl(2^{-1} (i+1)c(\beta )\lambda n^{1-\beta}2^{-n\beta}
\bigr)
\\
&&\hspace*{11pt}{}- S_\beta \bigl(2^{-1} ic(\beta )\lambda
n^{1-\beta}2^{-n\beta} \bigr) \leq2^{-n},\mbox{ some } i=0,\ldots, K_n \bigr)
\end{eqnarray*}
with $K_n \leq2 c(\beta)^{-1} n^{\beta-1}2^{n\beta}\log n$. Switching
to the
complements, and using once again (\ref{elowertail}) together with
the independence of the increments and self-similarity of the stable
subordinator, we conclude, after some straightforward calculus, that
for all $\lambda\geq\lambda(\beta)\in(0,\infty)$,
\[
\sum_{n=1}^\infty q ^{(2)}_n(
\lambda) \leq c(\beta)^{-1}\exp \bigl\{ -c(\beta)\lambda^{1/(1-\beta)}
\bigr\}.
\]
The resulting bound on the tail probability $P(K>\lambda)$ is
sufficient for the statement of the lemma.
\end{pf}

Recall that the only self-similar Gaussian process with stationary
increments is the Fractional Brownian motion (FBM), whose law is, apart from
the scale, uniquely determined by the self-similarity parameter $H\in
(0,1)$; see \citet{samorodnitskytaqqu1994}. This parameter of
self-similarity also determines the dependence properties of the
increment process of the FBM, the so-called fractional Gaussian noise,
with the case $H>1/2$ regarded as the long memory case. In contrast,
the self-similarity parameter almost never determines the dependence
properties of the increment processes of stable self-similar processes
with stationary increments; see
\citet{samorodnitsky2006LRD}. Therefore, it is interesting and
important to discuss the memory properties of the increment process
%
%
\begin{equation}
\label{eincrprocess} V_n^{(\alpha,\beta)} = Y_{\alpha,\beta}(n+1) -
Y_{\alpha,\beta}(n), \qquad n=0,1,2,\ldots.
\end{equation}
We refer the reader to \citet{rosinski1995} and
\citet{samorodnitsky2005a} for some of the notions used in the
statement of the following theorem.

%
\begin{theorem} \label{tincrprocess}
The stationary process
$ ( V_n^{(\alpha,\beta)}  )$ is generated by a conservative
null flow and is mixing.
\end{theorem}

\begin{pf}
Note that the increment process has the integral representation
\begin{eqnarray}
V_n^{(\alpha,\beta)} = \int_{\Omega^{\prime} \times[0,\infty)} \bigl(
M_{\beta} \bigl((n+1-x)_+,\omega^{\prime} \bigr) - M_{\beta}
\bigl((n-x)_+,\omega^{\prime} \bigr) \bigr) \,dZ_{\alpha,\beta}\bigl(
\omega^{\prime},x\bigr),\nonumber
\\\
\eqntext{n=0,1,2,\ldots.}
\end{eqnarray}
Since for every $x>0$, on a set of $P^\prime$ probability 1, by the
strong Markov property of the stable subordinator we have
\[
\limsup_{n\to\infty} M_{\beta} \bigl((n+1-x)_+ \bigr) -
M_{\beta} \bigl((n-x)_+ \bigr)>0,
\]
we see that
\[
\sum_{n=1}^\infty \bigl( M_{\beta}
\bigl((n+1-x)_+,\omega^{\prime
} \bigr) - M_{\beta} \bigl((n-x)_+,
\omega^{\prime} \bigr) \bigr)^\alpha= \infty, \qquad P^{\prime}
\times\nu\mbox{ a.e.}
\]
By Corollary 4.2 in \citet{rosinski1995}, we conclude that the
increment process is generated by a conservative flow.

It remains to prove that the increment process is mixing, since mixing
implies ergodicity which, in turn, implies that the
increment process is generated by a null flow; see
\citet{samorodnitsky2005a}. By Theorem 5
of \citet{rosinskizak1996}, it is enough to show that for every
$\epsilon> 0$,
\begin{eqnarray*}
&&\bigl(P^{\prime} \times\nu\bigr) \bigl\{ \bigl(\omega^{\prime}, x
\bigr)\dvtx  M_{\beta}\bigl((1-x)_+,\omega^{\prime}\bigr) > \epsilon,
\\
&&\hspace*{44pt} M_{\beta}\bigl((n+1-x)_+,\omega^{\prime}\bigr) - M_{\beta}
\bigl((n-x)_+,\omega^{\prime}\bigr) > \epsilon\bigr\}
\to0 \qquad\mbox{as } n \to\infty.
\end{eqnarray*}
However, an obvious upper bound on the expression in the left-hand
side is
\begin{eqnarray*}
&&\int_0^1 P^{\prime} \bigl(
M_{\beta}(n+1-x) - M_{\beta}(n-x) > \epsilon \bigr) (1-\beta)
x^{-\beta} \,dx
\\
&&\qquad =\int_0^1 P^{\prime}
\bigl(M_{\beta}\bigl((1-\delta_{n-x})_+\bigr) > \epsilon \bigr)
(1-\beta)x^{-\beta} \,dx,
\end{eqnarray*}
where for $r>0$, $\delta_{r}$ is a random variable, independent of the
Mittag--Leffler process, with the distribution given by
(\ref{eqDynkinLamperti}). Since $\delta_r$ converges weakly to
infinity as $r\to\infty$, by the dominated convergence theorem, the
above expression converges to zero as $n\to\infty$.
\end{pf}

%
\begin{remark} \label{rkbeta0}
Two extreme cases deserve mentioning.
A formal substitution of $\beta=0$ into (\ref{MTtransform}) leads
to a well-defined process $M_0(0)=0$ and $M_0(t)= E$, the same
standard exponential random variable for all $t>0$. This process is
no longer the inverse of a stable subordinator. It can, however, be
used in (\ref{eqMLSM}). It is elementary to see that the resulting
$\SaS$ process $Y_{\alpha,0}$ is, in fact, a $\SaS$ L\'evy motion.

On the other hand, a formal substitution of $\beta=1$ into
(\ref{MTtransform}) leads to the degenerate process $M_1(t)=t$ for
all $t\geq0$ [which can be viewed as the inverse of the degenerate
1-stable subordinator $S_1(t)=t$ for $t\geq0$]. Once again, this
process can be used in (\ref{eqMLSM}), if one interprets the measure
$\nu$ as the unit point mass at the origin. The resulting
$\SaS$ process $Y_{\alpha,1}$ is now the degenerate process
$Y_{\alpha,1}(t)=tY_{\alpha,1}(1)$ for all $t\geq0$, where
$Y_{\alpha,1}(1)$ is a $\SaS$ random variable.

Both limiting cases, $Y_{\alpha,0}$ and $Y_{\alpha,1}$, are processes
of a very different nature from the \mbox{$\beta$-}ML fractional $\SaS$ motion with $0<\beta<1$.
\end{remark}

\section{Some ergodic theory} \label{secergodic}

In this section, we present some elements of ergodic theory used in
this paper. The main reference for these notions is
\citet{aaronson1997}; see also \citet{zweimuller2009}.

Let $ (E,\mathcal{E}, \mu )$ be a $\sigma$-finite measure
space. We will often use the notation $A=B$ mod $\mu$ for $A,B\in
\mathcal{E}$ when $\mu(A\triangle B)=0$.

Let $T\dvtx   E \to E$ be a measurable map that preserves the
measure $\mu$. When the entire sequence $T, T^2,
T^3, \ldots$ of iterates of $T$ is involved, we will sometimes refer to
it as \textit{a~flow}. The map $T$ is called \textit{ergodic} if the only sets $A$ in
$\mathcal{E}$ for which $A=T^{-1}A$ mod $\mu$ are those for which
$\mu(A)=0$ or $\mu(A^c)=0$. The map $T$ is called \textit{conservative}
if
\[
\sum_{n=1}^\infty\mathbf{1}_A
\circ T^n=\infty\qquad\mbox{a.e. on }A
\]
for every $A\in\mathcal{E}$ with $\mu(A)>0$. If $T$ is ergodic, then
the qualification ``on $A$'' above is not needed.

\textit{The dual operator} $\widehat T$ is an operator $L^1(\mu)\to
L^1(\mu)$
defined by
\[
\widehat{T}f = \frac{d(\nu_f \circ T^{-1})}{d\mu}
\]
with $\nu_f$ a signed measure on $ (E,\mathcal{E}  )$ given
by $\nu_f(A) = \int_A f  \,d\mu$, $A\in\mathcal{E}$. The dual operator
satisfies the relation
%
%
\begin{equation}
\label{edualrel} \int_E (\widehat{T} f)\cdot g \,d\mu= \int
_E f\cdot(g \circ T) \,d\mu
\end{equation}
for $f\in L^1(\mu),   g\in L^\infty(\mu)$. For any nonnegative
measurable function $f$ on $E$, a similar definition gives a
nonnegative measurable function $\widehat{T} f$, and
(\ref{edualrel}) holds for any two nonnegative measurable functions
$f$ and $g$.

An ergodic conservative measure preserving map $T$ is called \textit{pointwise dual ergodic} if there is a sequence of positive constants
$a_n\to\infty$ such that
%
%
\begin{equation}
\label{epointwdualerg} \frac{1}{a_n}\sum_{k=1}^n
\widehat T^k f\to\int_E f \,d\mu\qquad
\mbox{a.e.}
\end{equation}
for every $f\in L^1(\mu)$. If the measure $\mu$ is infinite, pointwise
dual ergodicity rules out
invertibility of the map $T$; in fact, no factor of $T$ can be
invertible; see page~129 of \citet{aaronson1997}.

Sometimes the convergence of the type described in the definition
(\ref{epointwdualerg}) of pointwise dual ergodicity is uniform on
certain sets. Let $A\in{\mathcal E}$ be a set with
$0<\mu(A)<\infty$. We say that $A$ is \textit{a Darling--Kac set} for an
ergodic conservative measure preserving map $T$ if for some sequence
of positive constants $a_n\to\infty$,
%
%
\begin{equation}
\label{eDKset} \frac{1}{a_n}\sum_{k=1}^n
\widehat T^k \mathbf{1}_A\to\mu(A) \qquad
\mbox{uniformly, a.e. on }A
\end{equation}
[i.e., the convergence in (\ref{eDKset}) is uniform on a
measurable subset $B$ of $A$ with $\mu(B)=\mu(A)$]. By Proposition
3.7.5 of \citet{aaronson1997}, existence of a Darling--Kac set implies
pointwise dual ergodicity of $T$, so it is legitimate to use the same
sequence $(a_n)$ in (\ref{epointwdualerg}) and (\ref{eDKset}).

Given a set $A\in{\mathcal E}$, the map $\varphi\dvtx   E \to\mathbb{N}
\cup\{ \infty\}$ defined by $\varphi(x) = \inf\{n \geq1\dvtx   T^n x
\in
A \}$, $x\in E$ is called \textit{the first entrance time to $A$}. If $T$ is
conservative and ergodic (in addition to being measure preserving),
and $\mu(A)>0$, then $\varphi<\infty$
a.e. on $E$. It is natural to measure how often the set $A$ is
visited by the flow $(T^n)$ by \textit{the wandering rate} sequence
\[
w_n = \mu \Biggl(\bigcup_{k=0}^{n-1}
T^{-k} A \Biggr),\qquad n=1,2,\ldots.
\]
There are several alternative expressions for the wandering rate
sequence, the last two following from the fact that $T$ is measure
preserving:
%
%
\begin{eqnarray}\label{ewanderingalt}
w_n &=& \sum_{k=0}^{n-1}
\mu(A_k) = \sum_{k=0}^{n-1} \mu
\bigl( A \cap\{ \varphi> k \} \bigr)
\nonumber\\[-8pt]\\[-8pt]
& =& \sum_{k=1}^\infty
\min(k,n) \mu \bigl( A \cap\{ \varphi= k \} \bigr).\nonumber
\end{eqnarray}
Here, $A_0 = A$ and $A_k = A^c \cap\{\varphi= k \}$ for $k \geq
1$. If $\mu$ is an infinite measure, $T$ is
conservative and ergodic, and $0<\mu(A)<\infty$, then it follows from
(\ref{ewanderingalt}) that
%
%
\begin{equation}
\label{eqwanderingmu} w_n \sim\mu(\varphi< n) \qquad\mbox{as } n \to\infty.
\end{equation}

Let $T$ be a conservative ergodic measure preserving map. If a set $A$
is a Darling--Kac set, then there is a precise connection
between the return sequence $(w_n)$ and the normalizing sequence
$(a_n)$ in (\ref{eDKset}) [and hence, also in
(\ref{epointwdualerg})], assuming regular variation. Specifically,
if either $(w_n) \in RV_{1-\beta}$ or $(a_n) \in RV_{\beta}$ for some
$\beta\in[0,1]$, then
%
%
\begin{equation}
a_n \sim\frac{1}{\Gamma(2-\beta) \Gamma(1+\beta)} \frac{n}{w_n} \qquad\mbox{as } n \to
\infty. \label{eqprop387}
\end{equation}
Proposition 3.8.7 in \citet{aaronson1997} gives one direction of this
statement, but the argument is easily reversed.

We will also have an opportunity to use a variation of the notion of a
Darling--Kac set. Let $T$ be an ergodic conservative measure preserving
map. A set $A\in{\mathcal E}$ with $0<\mu(A)<\infty$
is said to be a uniform set for a nonnegative function \mbox{$g\in L^1(\mu)$}
if
%
%
\begin{equation}
\label{equniformset} \frac{1}{a_n}\sum_{k=1}^n
\widehat T^k g \to\int_E g \,d\mu \qquad
\mbox{uniformly, a.e. on }A.
\end{equation}
If $g=\mathbf{1}_A$, then a uniform set is just a Darling--Kac set.

\section{Central limit theorem associated with conservative null flows} \label{secCLT}
In this section, we state and discuss a functional central limit
theorem for stationary infinitely divisible processes generated by
certain conservative flows. Throughout, $T$ is an ergodic conservative
measure preserving map on an infinite $\sigma$-finite measure space
$ (E,\mathcal{E}, \mu )$, and $M$ a symmetric homogeneous
infinitely divisible random measure on $(E,\mathcal{E})$ with control
measure $\mu$ and local L\'evy measure $\rho$, satisfying the regular
variation with index $-\alpha$, $0<\alpha<2$ at infinity condition
(\ref{eregvarlevy}).
We will impose an extra assumption on the lower
tail of the local L\'evy measure: for some~$p_0<2$
%
%
\begin{equation}
\label{eqorginregularity} x^{p_0} \rho(x,\infty) \to0\qquad\mbox{as } x \to0.
\end{equation}

Let $f\dvtx E \to\mathbb{R}$ be a measurable function. We will assume that
$f$ is supported by a set of finite $\mu$-measure, and has the
following integrability properties:
%
%
\begin{equation}
\label{eqintegrabilitycond} f \in\cases{ L^{1\vee p}(\mu)\qquad\mbox{for some }
p>p_0, &\quad if $0<\alpha<1$,
\vspace*{2pt}\cr
L^\infty(\mu), &\quad
if $\alpha=1$,
\vspace*{2pt}\cr
L^2(\mu), &\quad if $1< \alpha<2$.}
\end{equation}
We will, further, assume that
%
%
\begin{equation}
\label{emeannonzero} \mu(f)=\int_E f(s) \mu(ds)\neq0.
\end{equation}

We consider a stochastic process $\mathbf{X}= ( X_1,X_2, \ldots)$
of the form (\ref{etheprocess})--(\ref{ethekernel}). The integral
is well defined under the condition
\[
\int_E\int_\mathbb{R}\min \bigl( 1,
x^2f_n(s)^2 \bigr) \rho(dx) \mu(ds)<\infty.
\]
It is not difficult to verify that this condition holds due
to the assumptions on the L\'evy measure $\rho$ and the integrability
conditions (\ref{eqintegrabilitycond}) on $f$. Therefore, the
process $\mathbf{X}$ is a well-defined infinitely divisible stochastic process. It is
automatically stationary. The L\'evy measure of each $X_n$ is given
by $\nu_{\mathrm{marg}}= (\rho\times\mu)\circ H^{-1}$, where $H\dvtx
\mathbb{R}\times E \to\mathbb{R}$ is given by $H(x,s) = xf(s)$. The
assumptions on
the L\'evy measure $\rho$ and the integrability conditions
(\ref{eqintegrabilitycond}) on $f$ imply that
\[
\nu_{\mathrm{marg}}(\lambda,\infty) \sim \biggl( \int_E
\bigl|f(s)\bigr|^\alpha \mu(ds) \biggr) \rho(\lambda,\infty)
\]
as $\lambda\to\infty$. It follows that the marginal tail of the
process itself is the same:
\[
P(X_n>\lambda)\sim \biggl( \int_E
\bigl|f(s)\bigr|^\alpha \mu(ds) \biggr) \rho(\lambda,\infty)
\]
as $\lambda\to\infty$; see \citet{rosinskisamorodnitsky1993}. In
particular, the margi\-nal distributions of the process $\mathbf{X}$ are
in the
domain of attraction of a $\SaS$ law; its memory is determined by the
operator $T$ through (\ref{ethekernel}).

We will assume that the operator $T$ has a Darling--Kac set $A$ [recall
(\ref{eDKset})], and that the normalizing sequence $(a_n)$ is
regularly varying with exponent $\beta\in(0,1)$. We will also
assume that the function $f$ is supported by $A$.
We will add an extra assumption on the set $A$. Recall the definition
of the set $A_n$ in (\ref{ewanderingalt}) as being the collection
of those points outside of $A$ that enter $A$ for the first time after
$n$ steps, $n=1,2,\ldots.$ We will assume that
there exists a measurable function $K\dvtx  E \to\mathbb{R}_+$ such that
%
%
\begin{equation}
\label{eqstrongDK} \frac{\widehat{T}{}^n \mathbf{1}_{A_n}}{\mu(A_n)} \to K\qquad\mbox {uniformly, a.e. on }A.
\end{equation}
This condition is an extension of the property shared by certain
operators $T$, the so-called Markov shifts [see Chapter~4 in
\citet{aaronson1997}], to a more general class of operators. See
Examples~\ref{exMC} and~\ref{exAFN} below.

Let $\rho^{\leftarrow}(y) = \inf \{ x\geq0\dvtx   \rho(x,\infty)
\leq
y \},
y>0$ be the left continuous inverse of the tail of the local
L\'evy measure. The regular variation of the tail implies that
$\rho^{\leftarrow}\in RV_{-1/\alpha}$ at zero. Define
%
%
\begin{equation}
\label{edefc} c_n = \Gamma(1+\beta) C_\alpha^{-1/\alpha}
a_n \rho ^{\leftarrow} (1/w_n),\qquad n=1,2,\ldots,
\end{equation}
where $C_\alpha$ is the $\alpha$-stable tail constant [see
\citet{samorodnitskytaqqu1994}],
$(a_n)$~is the normalizing sequence in the Darling--Kac property
(\ref{eDKset}) [or, equivalently, in the pointwise dual ergodicity
property (\ref{epointwdualerg})], and $(w_n)$ is the wandering
rate sequence for the set $A$ [related to the sequence $(a_n)$ via
(\ref{eqprop387})]. It follows immediately that
%
%
\begin{equation}
\label{ecnregvar} c_n \in RV_{\beta+ (1-\beta) / \alpha}.
\end{equation}
The sequence $(c_n)$ is the normalizing sequence in the functional
central limit theorem below. We will see that under the conditions of
that theorem we have the asymptotic relation
%
%
\begin{eqnarray}\label{eqasycn}
&& \rho ( c_n / a_n, \infty )\nonumber
\\[-3pt]
&&\qquad \sim C_\alpha \bigl(
C_{\alpha,\beta}/\Gamma(1+\beta) \bigr)^\alpha\bigl|\mu(f)\bigr|^\alpha
a_n^{\alpha} \Biggl( \int_E \Biggl| \sum
_{k=1}^n f \circ T^k
(x)\Biggr|^{\alpha} \mu(dx) \Biggr)^{-1}
\\
\eqntext{\mbox{as } n \to\infty}
\end{eqnarray}
with
%
%
\begin{equation}
\label{eCalphabeta} C_{\alpha,\beta} = \Gamma(1+\beta) \bigl( (1-\beta) B(1-\beta,
1+\alpha \beta) E \bigl(M_{\beta}(1)\bigr)^{\alpha}
\bigr)^{1/\alpha}.
\end{equation}
Here, $B$ is the standard beta function,
and $M_\beta$ the
Mittag--Leffler process defined in (\ref{eMLprocess}).
The following is our functional central limit theorem.

%
%
\begin{theorem} \label{tFCLT}
Let $T$ be an ergodic conservative
measure preserving map on an infinite $\sigma$-finite measure space
$ (E,\mathcal{E}, \mu )$, possessing a Darling--Kac set $A$
whose normalizing sequence $(a_n)$ is
regularly varying with exponent $\beta\in(0,1)$. Assume that
(\ref{eqstrongDK}) holds. Let $M$ be a symmetric homogeneous
infinitely divisible random measure on $(E,\mathcal{E})$ with control
measure $\mu$ and local L\'evy measure $\rho$, satisfying the regular
variation with index $-\alpha$, $0<\alpha<2$ at infinity condition~(\ref{eregvarlevy}). Assume further that
(\ref{eqorginregularity}) holds for some $p_0<2$.

Let $f$ be a measurable function supported by $A$ and satisfying
(\ref{eqintegrabilitycond}) and (\ref{emeannonzero}). If $1 <
\alpha< 2$, assume further
that either:
\begin{longlist}[(ii)]
\item[(i)] $A$ is a uniform set for $|f|$, or

\item[(ii)] $f$ is bounded.\vadjust{\goodbreak}
\end{longlist}

Then the stationary infinitely divisible stochastic process $\mathbf
{X}= ( X_1,X_2,\ldots)$
given by (\ref{etheprocess}) and (\ref{ethekernel}) satisfies
%
%
\begin{equation}
\label{emain} \frac{1}{c_n} \sum_{k=1}^{\lceil n\cdot\rceil}
X_k \Rightarrow \mu(f) Y_{\alpha, \beta} \qquad\mbox{in } D[0,\infty),
\end{equation}
where $(c_n)$ is defined by (\ref{edefc}), and $\{ Y_{\alpha,\beta}
\}$ is the \mbox{$\beta$-}Mittag--Leffler fractional $\SaS$ motion defined by (\ref{eqMLSM}).
\end{theorem}

%
\begin{remark}
The type of the limiting process obtained in Theorem
\ref{tFCLT} is an indication of the long memory in the process
$\mathbf{X}$. On the other hand, the Darling--Kac assumption
(\ref{eDKset}) and the duality relation (\ref{edualrel}) imply that
\begin{eqnarray*}
\frac{1}{a_n} \sum_{k=1}^n \mu
\bigl(A \cap T^{-k}A\bigr) &=& \frac{1}{a_n} \sum
_{k=1}^n \int_E
\mathbf{1}_A\cdot\mathbf{1}_A \circ T^k d
\mu
\\
&=& \int_A \frac{1}{a_n} \sum
_{k=1}^n \widehat T^k
\mathbf{1}_A \,d\mu \to\mu(A)^2 \in(0,\infty)
\end{eqnarray*}
as $n\to\infty$. Since $a_n=o(n)$, and $f$ is supported by $A$, we
see that for every $\epsilon> 0$,
\[
\frac{1}{n} \sum_{k=1}^n \mu
\bigl\{x \in E\dvtx  \bigl|f(x)\bigr| > \epsilon, \bigl|f \circ T^k(x)\bigr| > \epsilon \bigr\}
\leq\frac{1}{n} \sum_{k=1}^n \mu
\bigl(A \cap T^{-k}A\bigr) \to0
\]
and it follows immediately, for example, from Theorem 2 in
\citet{rosinskizak1997} that the process $\mathbf{X}$ is ergodic.

Under certain additional assumptions on the map $T$, one can check
that the process $\mathbf{X}$ is, in fact, mixing. We skip the
details. See,
however, Examples~\ref{exMC} and~\ref{exAFN} below.
\end{remark}

%
\begin{remark} \label{rkresultbetazero}
The statement of Theorem~\ref{tFCLT} makes sense in the limiting
cases $\beta=0$ and $\beta=1$ of Remark~\ref{rkbeta0} (in the case
$\beta=1$ the constant $C_{\alpha,1}$ needs to be interpreted as
$C_\alpha^{1/\alpha}$). Most of the argument in the proof of
Theorem~\ref{tFCLT} automatically works in these cases. The limiting
processes would then turn out to be, correspondingly, a $\SaS$ L\'evy
motion and the straight line process;
see Remark~\ref{rkbeta0}. This case $\beta=0$ corresponds to short
memory in the process $\mathbf{X}$, while the case $\beta=1$
corresponds to
extremely long memory.
\end{remark}

%
\begin{remark} \label{rkpositive}
When $0<\alpha<1$, the argument we will use in the proof of Theorem
\ref{tFCLT} can be used to establish a ``positive'' version of the
theorem. Specifically, assume now that the local L\'evy measure $\rho$
is concentrated on $(0,\infty)$, and that the function $f$ is
nonnegative. Then
%
%
\begin{equation}
\label{emainpos} \frac{1}{c_n} \sum_{k=1}^{\lceil n\cdot\rceil}
X_k \Rightarrow \mu(f) Y_{\alpha, \beta}^+ \qquad\mbox{in } D[0,
\infty),
\end{equation}
where $\{ Y_{\alpha,\beta}^+\}$ is a positive \mbox{$\beta$-}Mittag--Leffler
fractional $\alpha$-stable motion defined as in~(\ref{eqMLSM}), but
with $\SaS$ random measure $Z_{\alpha,\beta}$ replaced by a positive
$\alpha$-stable random measure with the same control measure.
\end{remark}

We finish this section with two examples of different situations where
Theorem~\ref{tFCLT} applies. The first example is close to the heart of
a probabilist.

%
\begin{example} \label{exMC}
Consider an irreducible null recurrent Markov chain with state space
$\mathbb{Z}$ and transition matrix $P=(p_{ij})$. Let $\{ \pi_j, j \in
\mathbb{Z} \}$ be the unique invariant measure of the Markov chain
that satisfies $\pi_0=1$. We define a $\sigma$-finite measure on
$(E,\mathcal{E}) = (\mathbb{Z}^{\mathbb{N}},
\mathcal{B}(\mathbb{Z}^{\mathbb{N}}) )$ by
\[
\mu(\cdot) = \sum_{i \in\mathbb{Z}} \pi_i
P_i(\cdot)
\]
with the usual notation of $P_i(\cdot)$ being the
probability law of the Markov chain starting in state $i \in
\mathbb{Z}$. Since $\sum_{j} \pi_j = \infty$, $\mu$ is an infinite
measure.

Let $T\dvtx   \mathbb{Z}^\mathbb{N}\to\mathbb{Z}^\mathbb{N}$ be the
left shift map
$T(x_0,x_1,\ldots) = (x_1,x_2, \ldots)$ for\break $\{x_k, k=0,1,\ldots\}\in
\mathbb{Z}^\mathbb{N}$. Obviously, $T$ preserves the measure $\mu$.
Since the
Markov chain is irreducible and null recurrent, the flow $\{T^n \}$ is
conservative and ergodic; see \citet{harrisrobbins1953}.

Consider the set $A =  \{ x \in\mathbb{Z}^{\mathbb{N}}\dvtx  x_0 = 0
\}$ and the corresponding first entrance time
$\varphi(x) = \min\{ n \geq1\dvtx    x_n=0 \}$,
$x\in\mathbb{Z}^\mathbb{N}$. Assume that
%
%
\begin{equation}
\label{esumprobregvar} \sum_{k=1}^n
P_0(\varphi\geq k) \in RV_{1-\beta}
\end{equation}
for some $\beta\in(0,1)$. Since $\mu(\varphi= k) = P_0(\varphi
\geq k)$ for $k\geq1$ [see Lemma 3.3 in
\citet{resnicksamorodnitskyxue2000}], we see that $\mu(\varphi\leq
n) \in RV_{1-\beta}$, and hence, by (\ref{eqwanderingmu}), the
wandering rates $(w_n)$ have the same property,
%
%
\begin{equation}
w_n \in RV_{1-\beta}. \label{eqMarkovwander}
\end{equation}

In this example,
\[
\widehat{T}{}^k \mathbf{1}_A (x) = P_0(x_k=0)
\qquad\mbox{constant for }x \in A;
\]
see Section~4.5 in \citet{aaronson1997}. In particular, the set $A$
is a Darling--Kac set, and by (\ref{eqMarkovwander}) and
(\ref{eqprop387}), we see that the corresponding normalizing
sequence $(a_n)$ is regularly varying with exponent $\beta$. Assumption
(\ref{eqstrongDK}) is easily seen to hold in this
example. Indeed, applying the explicit expression for the dual
operator given on page 156 in \citet{aaronson1997} to the function
\[
g(x_0,x_1,\ldots)= \mathbf{1} ( x_j\neq0,
j=0,\ldots, n-1, x_n=0 ),
\]
we see that
\[
\widehat{T}{}^n \mathbf{1}_{A_n} ( x_0,x_1,\ldots ) = \mathbf{1} (x_0=0 )\sum_{i_0\neq0}
\pi_{i_0}\sum_{i_1\neq0}p_{i_0
i_1}\cdots\sum_{i_{n-1}\neq0}p_{i_{n-2} i_{n-1}}p_{i_{n-1} 0}
\]
is constant on $A$ and vanishes outside of $A$. Therefore, the ratio
in (\ref{eqstrongDK}) is identically equal to 1 on $A$.

We conclude that Theorem~\ref{tFCLT} applies in this case if we
choose any
measurable function $f$ supported by $A$ and satisfying the conditions
of the theorem.

It is easy to see that the stationary infinitely divisible process
$\mathbf{X}$ in this
example is mixing. Indeed, by Theorem 5 of \citet{rosinskizak1996}
it is enough to check that
\[
\mu \bigl\{x\dvtx  \bigl|f(x)\bigr| > \epsilon, \bigl|f \circ T^n(x)\bigr| > \epsilon \bigr\}
\to0
\]
for every $\epsilon>0$. However, since $f$ vanishes outside of $A$,
null recurrence implies that as $n \to\infty$,
\[
\mu \bigl\{x\dvtx  \bigl|f(x)\bigr| > \epsilon, \bigl|f \circ T^n(x)\bigr| > \epsilon \bigr\}
\leq\mu\bigl(A \cap T^{-n}A\bigr) = P_0(x_n=0)
\to0.
\]
\end{example}

The next example is less familiar to probabilists, but is well known
to ergodic theorists.

%
\begin{example} \label{exAFN}
We start with a construction of the so-called \textit{AFN-system},
studied in, for example, \citet{zweimuller2000} and
\citet{thalerzweimuller2006}.
Let $E$ be the union of a finite family of disjoint bounded open
intervals in $\mathbb{R}$ and let $\mathcal{E}$ be the Borel $\sigma$-field
on $E$. Let $\lambda$ be the one-dimensional Lebesgue measure.

Let $\xi$ be a (possibly, infinite) collection of nonempty disjoint open
subintervals (of the intervals in $E$) such that $\lambda ( E
\setminus
\bigcup_{ Z\in\xi}Z ) = 0$. Let $T\dvtx E \to E$ be a map that is twice
differentiable on (each interval of) $E$. We assume that $T$ is
strictly monotone on each $Z \in\xi$.

Map $T$ is further assumed to satisfy the following three
conditions, (A), (F) and~(N), (giving rise to the name
AFN-system).
\begin{longlist}[(B)]
\item[(A)]  \textit{Adler's condition}:
\[
T^{\prime\prime} / \bigl(T^{\prime}\bigr)^2 \mbox{ is
bounded on } \bigcup_{ Z\in\xi}Z.
\]

\item[(F)]  \textit{Finite image condition}:
\[
\mbox{the collection } T \xi= \{ TZ\dvtx  Z \in\xi\}\mbox{ is finite}.
\]
\item[(N)] \textit{A possibility of nonuniform expansion}:
there exists a finite subset $\zeta\subseteq\xi$ such that each $Z
\in\zeta$ has \textit{an indifferent fixed-point} $x_Z$ as one of its end
points. That is,
\[
\lim_{x \to x_Z, x \in Z} Tx = x_Z \quad\mbox{and}\quad\lim
_{x \to
x_Z, x \in Z} T^{\prime} x = 1.
\]
Moreover, we suppose, for each $Z \in\zeta$,
\[
\mbox{either } T^{\prime} \mbox{ decreases on } (-\infty,
x_Z) \cap Z\quad\mbox{or}\quad T^\prime\mbox{ increases on }
(x_Z, \infty) \cap Z,
\]
depending on whether $x_Z$ is the left endpoint or the right endpoint
of $Z$.
Finally, we assume that $T$ is uniformly expanding away from $\{x_Z\dvtx
Z \in\zeta\}$, that is, for each $\epsilon> 0$, there is
$\rho(\epsilon) > 1$ such that
\[
\bigl|T^{\prime}\bigr| \geq\rho(\epsilon)\qquad\mbox{on } E \bigm\backslash \bigcup
_{Z
\in\zeta} \bigl( (x_Z - \epsilon,
x_Z + \epsilon ) \cap Z \bigr).
\]
\end{longlist}

If the conditions (A), (F), and (N) are satisfied, the triplet
$(E,T,\xi)$ is called an AFN-system, and the map $T$ is called an
AFN-map. If $T$ is also conservative and ergodic with respect to
$\lambda$,
and the collection $\zeta$ is nonempty, then the AFN-map $T$ is said
to be \textit{basic}; we will assume this property in the
sequel. Finally, we will assume that $T$ admits \textit{nice expansions}
at the indifferent fixed points. That is, for every $Z\in\zeta$ there
is $0 < \beta_Z < 1$ such that
%
%
\begin{equation}
\label{eniceexp} \qquad\quad Tx = x + a_Z |x-x_Z|^{1/\beta_Z + 1} +
o \bigl(|x-x_Z|^{1/\beta_Z +
1} \bigr)\qquad\mbox{as } x \to
x_Z\mbox{ in } Z
\end{equation}
for some $a_Z \neq0$.

It is shown in \citet{zweimuller2000} that every basic AFN-map has an
infinite invariant measure $\mu\ll\lambda$ with the density given by
$d\mu/ d \lambda(x) = h_0(x) G(x)$, \mbox{$x\in E$,} where
\[
G(x) = %
\cases{\displaystyle (x-x_Z) \bigl(x-(T|_Z)^{-1}(x)
\bigr)^{-1}, &\quad if $x \in Z \in\zeta$,
\vspace*{3pt}\cr
1, &\quad if $\displaystyle x \in E \bigm\backslash\bigcup_{Z\in\zeta}Z$}
\]
and $h_0$ is a function of bounded variation bounded away from both 0
and infinity. We view $T$
as a conservative ergodic measure-preserving map on the infinite
measure space $(E,\mathcal{E},\mu)$.

An example of a basic AFN-map is Boole's transformation placed on
$E=(0,1/2)\cup(1/2,1)$, defined by
\begin{eqnarray*}
T(x) &=& \frac{x(1-x)}{1-x-x^2},\qquad x\in(0,1/2),
\\
T(x) &=& 1-T(1-x),\qquad x\in(1/2,1).
\end{eqnarray*}
It admits nice expansions at the indifferent fixed points $x_Z=0$ and
$x_Z=1$ with $\beta_Z=1/2$ in both cases. The invariant measure $\mu$
satisfies
\[
\frac{d\mu}{d\lambda}(x) = \frac{1}{x^2}+ \frac{1}{(1-x)^2},\qquad x\in E.
\]
See \citet{thaler2001}.

Let $T$ be a basic AFN-map. We put
\[
A=E \bigm\backslash\bigcup_{Z
\in\zeta} \bigl(
(x_Z - \epsilon, x_Z + \epsilon ) \cap Z \bigr)
\]
for some $\epsilon> 0$ small enough so that the set $A$ is
nonempty. Since $\lambda(\partial
A)=0$ and $A$ is bounded away from the indifferent fixed points
$\{ x_Z\dvtx  Z \in\zeta\}$, it follows from Corollary 3 of
\citet{zweimuller2000} that $A$ is a Darling--Kac set. Moreover, the
corresponding normalizing sequence $(a_n)$ is regularly varying with
exponent $\beta= \min_{Z \in\zeta} \beta_Z$ in the notation of
(\ref{eniceexp}); see Theorems 3 and 4 in
\citet{zweimuller2000}. The assumption (\ref{eqstrongDK}) also holds;
see (2.6) in \citet{thalerzweimuller2006}.

Once again, Theorem~\ref{tFCLT} applies if we choose any measurable
function $f$ supported by $A$ and satisfying the
conditions of Theorem~\ref{tFCLT}. Note that, by Theorem~9 in
\citet{zweimuller2000}, Riemann integrability of $|f|$ on $A$ suffices
for the uniformity of the set $A$ for $|f|$.

The stationary infinitely divisible process $\mathbf{X}$ in this
example is also
mixing. Indeed, the basic AFN-map $T$ is \textit{exact}, that is, the
$\sigma$-field $\bigcap_{n=1}^\infty T^{-n}{\mathcal{E}}$ is trivial; see,
for example, page 1522 in \citet{zweimuller2000}. The exactness of
$T$ implies
that
\[
\mu\bigl(A \cap T^{-n}A\bigr) = \int_A
\widehat{T}{}^n \mathbf{1}_{A} \,d\mu \to0
\]
as $n\to\infty$; see page 12 in \citet{thaler2001}. Now mixing of the
process $\mathbf{X}$ follows from the fact that $f$ is supported by
$A$, as
in Example~\ref{exMC}.
\end{example}

\section{Distributional results in ergodic theory} \label{secdistrresults}

In this section, we prove two distributional ergodic-theoretical
results that will be used in the proof of Theorem~\ref{tFCLT}. These
results may be of interest on their own as well. We call our first
result a generalized Darling--Kac theorem, because the first result of
this type was proved in \citet{darlingkac1957} as a distributional
limit theorem for the occupation times of Markov processes and chains
under a certain uniformity assumption on the transition law. The
limiting law is the Mittag--Leffler distribution described in
(\ref{MTtransform}). Under the same setup and assumptions,
\citet{bingham1971} extended the result to weak convergence in the
space $D[0,\infty)$ endowed with the Skorohod $J_1$ topology, and the
limiting process is the Mittag--Leffler process defined in
(\ref{eMLprocess}).

The result of \citet{darlingkac1957} was put into ergodic-theoretic
context by \citet{aaronson1981} who established the one-dimensional
convergence for abstract conservative infinite measure preserving maps
under the assumption of pointwise dual ergodicity, that is, dispensing
with a condition of uniformity. Furthermore, Aaronson proves
convergence in a \textit{strong distributional sense}, a stronger mode
of convergence than weak convergence. The same strong distributional
convergence was established later in \citet{thalerzweimuller2006},
with the assumption of pointwise dual ergodicity replaced by an
averaged version of (\ref{eqstrongDK}). The latter assumption was further
weakened in \citet{zweimuller2007}. Our result, Theorem
\ref{tgenDK} below, extends Aaronson's result to the space
$D[0,\infty)$, under the assumption of pointwise dual ergodicity.

We start with defining strong distributional convergence. Let $Y$ be a
separable metric space, equipped with its Borel $\sigma$-field. Let
$ ( \Omega_1, {\mathcal F}_1, m )$ be a measure space and
$ ( \Omega_2, {\mathcal F}_2, P_2 )$ a probability space. We
say that a sequence of measurable maps $R_n\dvtx   \Omega_1\to Y$,
$n=1,2,\ldots$ converges strongly in distribution to a measurable map
$R\dvtx   \Omega_2\to Y$ if $P_1\circ R_n^{-1}\Rightarrow P_2\circ R^{-1}$ in
$Y$ for any probability measure $P_1\ll m$ on $ ( \Omega_1,
{\mathcal F}_1 )$. That is,
\[
\int_{\Omega_1} g(R_n) \,dP_1 \to\int
_{\Omega_2} g(R) \,dP_2
\]
for any such $P_1$ and a bounded continuous function $g$ on $Y$. We
will use the notation $R_n \stackrel{\mathcal{L}(m)}{\Rightarrow}
R$ when strong distributional convergence takes place.

%
\begin{theorem}[(Generalized Darling--Kac theorem)]\label{tgenDK}
Let $T$ be an ergodic conservative
measure preserving map on an infinite $\sigma$-finite measure space
$ (E,\mathcal{E}, \mu )$. Assume that $T$ is pointwise dual
ergodic with a normalizing sequence $(a_n)$ that is
regularly varying with exponent $\beta\in(0,1)$. Let $f\in
L^1(\mu)$ be such that $\mu(f)\neq0$, and denote $S_n(f) =
\sum_{k=1}^n f \circ T^k$, $n=1,2,\ldots.$ Then
%
%
\begin{equation}
\frac{1}{a_n} S_{\lceil n \cdot\rceil} (f) \stackrel{\mathcal{L}(\mu)} {
\Rightarrow} \mu(f) \Gamma(1+\beta) M_{\beta}(\cdot) \qquad\mbox{in } D[0,
\infty), \label{eqgoalLemma1}
\end{equation}
where $M_\beta$ is the Mittag--Leffler process, and $D[0,\infty)$ is
equipped with the $J_1$ topology.
\end{theorem}

\begin{pf}
It is shown in Corollary 3 of \citet{zweimuller2007a} that proving weak
convergence in (\ref{eqgoalLemma1}) for one fixed probability
measure on $ (E,\mathcal{E} )$, that is absolutely continuous
with respect to $\mu$, already guarantees the full strong
distributional convergence. We choose and fix an arbitrary set $A\in
\mathcal{E}$ with $0<\mu(A)<\infty$, and prove weak
convergence in (\ref{eqgoalLemma1}) with respect to $\mu_A(\cdot)
=\mu(\cdot\cap A) / \mu(A)$.

It turns out that we only need to consider one particular function
$f=\mathbf{1}_A$ and to establish the appropriate finite-dimensional
convergence, that is, to show that
%
%
\begin{equation}
\biggl( \frac{1}{a_n} S_{\lceil n t_i \rceil} (\mathbf{1}_A)
\biggr)_{i=1}^k \Rightarrow \bigl( \mu(A) \Gamma(1+\beta)
M_{\beta}(t_i) \bigr)_{i=1}^k \qquad
\mbox{in }\mathbb{R}^k \label{eqlonggoal}
\end{equation}
for all $k \geq1$, $0 \leq t_1 < \cdots< t_k$, when the law of the
random vector in the left-hand side is computed with respect to
$\mu_A$.

Indeed, suppose that (\ref{eqlonggoal}) holds. By Hopf's ergodic
theorem [also sometimes called a ratio ergodic theorem; see Theorem
2.2.5 in \citet{aaronson1997}], the finite-dimensional convergence
immediately extends to the corresponding finite-dimensional
convergence with any function $f\in L^1(\mu)$ such that $\mu(f)\neq
0$. Next, write $f=f_+-f_-$, the
difference of the positive and negative parts. Since the process
$ ( a_n^{-1}S_{\lceil nt \rceil}(f_+),   t\geq0 )$ has,
for each
$n$, nondecreasing sample paths, Theorem~3 in \citet{bingham1971}
tells us that the convergence of the
finite-dimensional distributions, and the continuity in probability of
the limiting Mittag--Leffler process already imply weak convergence,
hence tightness, of this sequence of processes. Similarly, the
sequence of the processes $ ( a_n^{-1}S_{\lceil nt \rceil}(f_-),
t\geq
0 )$, $n=1,2,\ldots$ is tight as well. Since both converge to a
continuous limit, their sum, $ ( a_n^{-1}S_{\lceil nt \rceil}(f),
t\geq
0 )$, $n=1,2,\ldots,$ is tight as well, because in this case the
uniform modulus of continuity can be used instead of the $J_1$ modulus
of continuity; see, for example, \citet{billingsley1999}.

This will give us the required weak convergence, and hence, complete the
proof of the theorem.

It remains to show (\ref{eqlonggoal}). We will use a strategy similar
to the one used in \citet{bingham1971}. We start with defining a
continuous version of the process $ ( S_{\lceil n t \rceil}
(\mathbf{1}_A),
t\geq0 )$ given by the linear interpolation
%
%
\begin{eqnarray}\label{einterpprocess}
\widetilde S_n(t) &=& \bigl( (i+1)-nt \bigr)S_i(
\mathbf{1}_A) + (nt-i)S_{i+1}(\mathbf{1} _A)
\nonumber\\[-8pt]\\[-8pt]
\eqntext{\displaystyle\mbox{if } \frac{i}{n}\leq t\leq\frac{i+1}{n}, i=0,1,2,\ldots.}
\end{eqnarray}
With the implicit argument $x\in E$ viewed as random (with the law
$\mu_A$), each $\widetilde S_n$ defines a random Radon measure on
$[0,\infty)$. Therefore, for any $k\geq1$ the $k$-tuple product
$\widetilde S_n^k = \widetilde S_n\times\cdots\times\widetilde S_n$ is a random
Radon measure on $[0,\infty)^k$. By Fubini's theorem,
\[
\tilde m_n^{(k)}(B) = \int_A \widetilde S_n^k(B) (x) \mu_A(dx),\qquad B\subseteq[0,
\infty)^k, \mbox{ Borel,}
\]
is a Radon measure on $[0,\infty)^k$. We define, similarly,
$S_n$, $S_n^k$ and $m_n^{(k)}$, starting with $S_n(t) = S_{\lceil n t
\rceil} (\mathbf{1}_A),   t\geq0$. Finally, we perform the same
operation on the limiting process and define $M_{\beta,A}$ by
$\mu(A)\Gamma(1+\beta)M_\beta$, and then construct $M_{\beta,A}^k$ and
$m_{\beta,A}^{(k)}=E
M_{\beta,A}^k$.

Note that $\tilde m_n^{(k)}$ is
absolutely continuous with respect to the $k$-dimensional Lebesgue
measure, and
\begin{eqnarray}
\frac{d^k \tilde m_n^{(k)}}{dt_1\cdots\, dt_k} &=& n^k \int_A \prod
_{j=1}^k \mathbf{1}_A \circ
T^{i_j}(x) \mu_A(dx)\nonumber
\\
\eqntext{\displaystyle\mbox{on } \frac{i_j}{n}
\leq t_j<\frac{i_j+1}{n}, i_j=0,1,\ldots, j=1,\ldots, k.}
\end{eqnarray}

We will prove that for all $k \geq1$,
$\theta_1, \ldots, \theta_k \geq0$,
%
%
\begin{eqnarray}\label{eqLSTconv}
&& \frac{1}{a_n^k} \int_0^{\infty}
\cdots\int_0^{\infty} e^{-\sum_{j=1}^k \theta_j t_j} \tilde
m_n^{(k)}(dt_1\cdots \, dt_k)
\nonumber\\[-8pt]\\[-8pt]
&&\qquad \to\int
_0^{\infty}   \cdots\int
_0^{\infty} e^{-\sum_{j=1}^k \theta_j t_j} m_{\beta,A}^{(k)}(dt_1\cdots \, dt_k) \nonumber
\end{eqnarray}
as $n\to\infty$. We claim that this will suffice for
(\ref{eqlonggoal}).

Indeed, suppose that (\ref{eqLSTconv}) holds. Convergence of the
joint Laplace transforms implies that
\[
a_n^{-k}\tilde m_n^{(k)} \stackrel{v}
{\rightarrow} m_{\beta,A}^{(k)}
\]
(vaguely) in $[0,\infty)^k$. Since the rectangles are, clearly,
compact continuity sets with respect to the limiting measure
$m_{\beta,A}^{(k)}$, we conclude that for every $k=1,2,\ldots$ and
$t_j\geq0,   j=1,\ldots, k$, we have
\begin{eqnarray*}
&& \int_A \prod_{j=1}^k
a_n^{-1}\widetilde S_n(t_j) (x)
\mu_A(dx)
\\
&&\qquad
= a_n^{-k}\tilde
m_n^{(k)} \Biggl( \prod_{j=1}^k[0,t_j]
\Biggr)\to m_{\beta,A}^{(k)} \Biggl( \prod_{j=1}^k[0,t_j]
\Biggr)
\\
&&\qquad
= E \Biggl[ \prod_{j=1}^k \mu(A)
\Gamma(1+\beta) M_\beta(t_j) \Biggr]
\end{eqnarray*}
as $n\to\infty$. Since for every fixed $\varepsilon>0$ and
$n>1/\varepsilon$,
\[
\widetilde S_n(t)\leq S_n(t)\leq\widetilde S_n(t+
\varepsilon)
\]
for each $t\geq0$, we conclude by monotonicity and continuity of the
Mittag--Leffler process that
%
%
\begin{equation}
\label{eprodconvtilde} \int_A \prod
_{j=1}^k a_n^{-1}S_n(t_j)
\mu_A(dx) \to E \Biggl[ \prod_{j=1}^k
\mu(A) \Gamma(1+\beta) M_\beta (t_j) \Biggr].
\end{equation}
We claim that (\ref{eprodconvtilde}) implies
(\ref{eqlonggoal}). By taking linear combinations with nonnegative
weights, we see that it is enough to show that the distribution of
such a linear combination,
\[
\sum_{j=1}^k \theta_jM_\beta(t_j),\qquad
\theta_j>0, j=1,\ldots, k,
\]
is determined by its moments, and by the Carleman sufficient condition
it is enough to check that
\[
\sum_{m=1}^\infty \biggl(\frac{1}{E ( \sum_{j=1}^k
\theta_jM_\beta(t_j) )^m}
\biggr)^{1/(2m)} = \infty.
\]
A simple monotonicity and scaling argument shows that it is sufficient
to verify only that
%
%
\begin{equation}
\label{ecarleman} \sum_{m=1}^\infty \biggl(
\frac{1}{E (
M_\beta(1) )^m} \biggr)^{1/(2m)} = \infty.
\end{equation}
However, the moments of $M_\beta(1)$ can be read off (\ref
{MTtransform}), and Stirling's formula together with elementary
algebra imply (\ref{ecarleman}). Hence, (\ref{eqlonggoal}) follows.

It follows that we need to
prove (\ref{eqLSTconv}). Taking into account the form of the density
of $\tilde m_n^{(k)}$ with respect to the $k$-dimensional Lebesgue
measure, we can write the left-hand side of (\ref{eqLSTconv}) as
\[
\sum_{\pi} F_{n,A} (\theta_{\pi(1)}\cdots\theta_{\pi(k)}),
\]
where
\begin{eqnarray*}
&& F_{n,A} (\theta_1 \cdots\theta_k)
\\
&&\qquad = \biggl(
\frac{n}{a_n} \biggr)^k \idotsint_{0<t_1 < \cdots< t_k}
e^{-\sum_{j=1}^k
\theta_j t_j} \mu_A \Biggl( \bigcap
_{j=1}^k T^{- \lceil nt_j \rceil} A \Biggr)
\,dt_1\cdots \, dt_k
\end{eqnarray*}
and $\pi$ runs through the permutations of the sets $\{1,\ldots,k
\}$. To establish (\ref{eqLSTconv}), it is enough to verify that
%
%
\begin{eqnarray}\label{eqFnrelation}
&& F_{n,A}(\theta_1 \cdots\theta_k)
\nonumber\\[-8pt]\\[-8pt]
&&\qquad  \to \bigl(
\mu(A) \Gamma(1+\beta ) \bigr)^k \bigl((\theta_1 + \cdots+
\theta_k) (\theta_2 + \cdots+ \theta_k) \cdots
\theta_k \bigr)^{-\beta} \nonumber
\end{eqnarray}
as $n \to\infty$, because Lemma 3 in \citet{bingham1971} shows that
summing up the expression in the right-hand side of
(\ref{eqFnrelation}) over all possible
permutations $(\theta_{\pi(1)} \cdots\theta_{\pi(k)})$ produces the
expression in the right-hand side of (\ref{eqLSTconv}).

Given $0<\varepsilon<1$, we use repeatedly pointwise dual ergodicity and
Egorov's theorem to construct a nested sequence of measurable subsets
of $E$, with $A_0=A$, and for $i=0,1,\ldots,$ $A_{i+1}\subseteq A_i$,
and $\mu(A_{i+1})\geq(1-\varepsilon)\mu(A_i)$, while
%
%
\begin{equation}
\label{eapproxDK} \frac{1}{a_n}\sum_{k=1}^n
\widehat T^k \mathbf{1}_{A_i}\to\mu(A_i)
\qquad\mbox{uniformly on }A_{i+1}.
\end{equation}
It is elementary to see that with $v_1=\theta_1 + \theta_2+\cdots+
\theta_k$, $v_2=\theta_2 + \cdots+ \theta_k, \ldots, v_k=\theta_k$,
%
%
\begin{eqnarray}\label{ediscrF}
&& F_{n,A} (\theta_1 \cdots
\theta_k)\nonumber
\\
&&\qquad \sim\frac{1}{a_n^k} \sum_{m_1=0}^\infty
\cdots\sum_{m_k=0}^\infty e^{-n^{-1}\sum_{j=1}^k v_j m_j}
\mu_A \Biggl( \bigcap_{j=1}^k
T^{- (m_1+\cdots+m_j) } A \Biggr)\nonumber
\\
&&\qquad = \frac{1}{a_n^k} \bigintsss_A \Biggl[
\Biggl( \sum_{m_1=0}^\infty \widehat
T^{m_1} \mathbf{1}_{A} e^{- v_1m_1/n} \Biggr)
\\
&&\hspace*{63pt}{}\times \prod
_{j=2}^k \Biggl( \sum
_{m_j=0}^\infty\mathbf{1}_A\circ
T^{m_2+\cdots+m_j} e^{-
v_jm_j/n} \Biggr) \Biggr] \,d\mu_A\nonumber
\\
&&\qquad \geq\frac{1}{a_n^k} \int_{A_1} ( \cdots ),\nonumber
\end{eqnarray}
where the equality is due to the duality relation (\ref{edualrel}).
Note that by (\ref{eapproxDK}) with $i=0$,
%
%
\begin{eqnarray}\label{euseDK}
&& \sum_{m_1=0}^{\infty}
\widehat{T}{}^{m_1} \mathbf{1}_A e^{- v_1m_1/n}
\nonumber\\[-8pt]\\[-8pt]
&&\qquad =
\bigl(1-e^{-v_1/n} \bigr) \sum_{i=0}^{\infty}
\Biggl( \sum_{m_1=0}^i
\widehat{T}{}^{m_1} \mathbf{1}_{A_0} \Biggr) e^{- v_1 i /n}
\sim\frac{\mu(A_0)v_1}{n} \sum_{i=0}^{\infty}
a_i e^{- v_1 i/n}\hspace*{-30pt}
\nonumber
\end{eqnarray}
uniformly on $A_1$ as $n \to\infty$. Therefore,
\begin{eqnarray*}
&& F_{n,A} (\theta_1 \cdots\theta_k)
\\[-1pt]
&&\qquad \geq  \bigl(
1-o(1) \bigr) \frac{1}{a_n^k} \frac{\mu(A_0)v_1}{n}
\\[-1pt]
&&\quad\qquad{}\times \sum
_{i=0}^{\infty} a_i e^{-
v_1i/n} \bigintsss_{A_1}\prod_{j=2}^k
\Biggl( \sum_{m_j=0}^\infty
\mathbf{1}_A\circ T^{m_2+\cdots+m_j} e^{-
v_jm_j/n} \Biggr) \,d
\mu_A
\\[-1pt]
&&\qquad = \bigl( 1-o(1) \bigr) \frac{1}{a_n^k} \frac{\mu(A_0)v_1}{n}
\\[-1pt]
&&\quad\qquad{}\times  \sum
_{i=0}^{\infty} a_i e^{-
v_1i/n}
\bigintsss_{A} \Biggl[ \Biggl( \sum
_{m_2=0}^\infty \widehat T^{m_2}
\mathbf{1}_{A_1} e^{- v_2m_2/n} \Biggr)
\\[-1pt]
&&\hspace*{120pt}{}\times \prod
_{j=3}^k \Biggl( \sum_{m_j=0}^\infty
\mathbf{1}_A\circ T^{m_3+\cdots+m_j} e^{-
v_jm_j/n} \Biggr) \Biggr]
\,d\mu_A
\\[-1pt]
&&\qquad \geq \bigl( 1-o(1) \bigr) \frac{1}{a_n^k} \frac{\mu(A_0)v_1}{n} \sum
_{i=0}^{\infty} a_i e^{-
v_1i/n} \int
_{A_2} ( \cdots ).
\end{eqnarray*}
Using now repeatedly (\ref{eapproxDK}) with larger and larger $i$,
together with the same argument as in (\ref{euseDK}), we conclude
that
\begin{eqnarray*}
&& F_{n,A} (\theta_1 \cdots\theta_k)
\\[-1pt]
&&\qquad \geq  \bigl(
1-o(1) \bigr) \frac{1}{a_n^k} \frac{\mu(A_0)\mu(A_1)v_1v_2}{n^2}
\sum
_{i_1=0}^{\infty
} a_{i_1} e^{-
v_{1}i_1/n}\sum
_{i_2=0}^{\infty} a_{i_2}
e^{-
v_2i_2/n}
\\[-1pt]
&&\quad\qquad{}\times
\bigintsss_{A_2}\prod_{j=3}^k
\Biggl( \sum_{m_j=0}^\infty
\mathbf{1}_A\circ T^{m_3+\cdots+m_j} e^{-
v_jm_j/n} \Biggr) \,d
\mu_A
\\[-1pt]
&&\qquad \geq \cdots\geq \bigl( 1-o(1) \bigr) \frac{1}{a_n^k} \frac{\prod_{j=0}^{k-1}\mu(A_j)v_{j+1}}{n^k} \prod
_{j=1}^k \Biggl( \sum
_{i=0}^{\infty} a_i e^{- v_ji/n} \Biggr)
\frac{\mu(A_k)}{\mu(A)}
\\[-1pt]
&&\qquad \geq \bigl( 1-o(1) \bigr) (1-\varepsilon)^{k(k+1)/2} \biggl(
\frac{\mu(A)}{na_n} \biggr)^k (v_1\cdots v_k)
\prod_{j=1}^k \Biggl( \sum
_{i=0}^{\infty} a_i e^{-
v_ji/n} \Biggr).
\end{eqnarray*}
Extending the sequence $(a_n)$ into a piece-wise constant regular
varying function of real variable $ ( a(x),   x>0 )$ and
using Karamata's Tauberian theorem [see, e.g., Section~3.6 in
\citet{aaronson1997}], we conclude that for every $j=1,\ldots, k$,
\[
\sum_{i=0}^{\infty} a_i
e^{- v_ji/n} \sim\Gamma(1+\beta) \frac{n}{v_j} a(n/v_j),
\qquad n\to\infty.
\]
It follows that
\begin{eqnarray*}
F_{n,A} (\theta_1 \cdots\theta_k) &\geq& \bigl(
1-o(1) \bigr) (1-\varepsilon)^{k(k+1)/2} \bigl( \mu(A) \Gamma(1+\beta)
\bigr)^k \prod_{j=1}^k
\frac{a(n/v_j)}{a_n}
\\
&\to&(1-\varepsilon)^{k(k+1)/2} \bigl( \mu(A) \Gamma(1+\beta)
\bigr)^k \prod_{j=1}^k
v_j^{-\beta}
\end{eqnarray*}
by the regular variation. Since this is true for every $0<\varepsilon
<1$, we
have obtained the lower bound
%
%
\begin{eqnarray}\label{elowerbound}
&& \liminf_{n\to\infty} F_{n,A}(
\theta_1 \cdots\theta_k)
\nonumber\\[-8pt]\\[-8pt]
&&\qquad \geq \bigl(\mu(A) \Gamma(1+\beta)
\bigr)^k \bigl((\theta_1 + \cdots+ \theta_k) (
\theta_2 + \cdots+ \theta_k) \cdots \theta_k
\bigr)^{-\beta}.\nonumber
\end{eqnarray}

The lower bound (\ref{elowerbound}) is valid for any measurable set
$A$ with $0<\mu(A)<\infty$. We will now show that for any $k\geq1$
and $0<\theta<1$ there is a measurable set $A_{k,\theta}\subseteq A$ such
that
%
%
\begin{equation}
\label{ealmostA} \mu ( A_{k,\theta} ) \geq(1-\theta)\mu(A)
\end{equation}
and such that
%
%
\begin{eqnarray}\label{eupperbound}
&& \limsup_{n\to\infty} F_{n,A_{k,\theta}}(
\theta_1 \cdots\theta_k)
\nonumber\\[-8pt]\\[-8pt]
&&\qquad \leq \bigl(\mu(A_{k,\theta})
\Gamma(1+\beta) \bigr)^k \bigl((\theta_1 + \cdots+
\theta_k) (\theta_2 + \cdots+ \theta_k) \cdots\theta_k \bigr)^{-\beta}.\nonumber
\end{eqnarray}
We know that (\ref{elowerbound}) and (\ref{eupperbound}) together
imply (\ref{eqFnrelation}), hence that (\ref{eqlonggoal}) holds for
the set $A_{k,\theta}$. We claim that this implies that
(\ref{eqlonggoal}) for every measurable $A$ with $0<\mu(A)<\infty$.

Indeed, suppose that, to the contrary, (\ref{eqlonggoal}) fails for
some measurable $A$ with $0<\mu(A)<\infty$, some $k\geq1$ and
some $0<t_1<\cdots<t_k$. By the one-dimensional result of
\citet{aaronson1981}, the $k$ components in the left hand side of
(\ref{eqlonggoal}), individually, converge weakly. Therefore, the
sequence of the laws of the $k$-dimensional vectors in the left-hand
side of (\ref{eqlonggoal}) is tight, and so there is a sequence of
integers $n_l\uparrow\infty$ and a random vector $(Y_1,\ldots, Y_k)$
with
%
%
\begin{equation}
\label{ewronglimit} (Y_1,\ldots, Y_k)\stackrel{d} {\neq}
\mu(A) \Gamma(1+\beta) \bigl( M_{\beta}(t_1)\cdots M_{\beta}(t_k) \bigr),
\end{equation}
such that
%
%
\begin{equation}
\label{epresumeconv} \frac{1}{a_{n_l}} \bigl( S_{\lceil n_l t_1 \rceil} (
\mathbf{1}_A), \ldots, S_{\lceil n_l t_k \rceil} (\mathbf{1}_A)
\bigr)\Rightarrow (Y_1,\ldots, Y_k),
\end{equation}
when the law of the
random vector in the left-hand side is computed with respect to~$\mu_A$. It follows from (\ref{ewronglimit}) that there is a Borel set
$B\subset\mathbb{R}^k$ such that, for each $b>0$, $bB$ is a
continuity set
for both $(Y_1,\ldots, Y_k)$ and $ \mu(A) \Gamma(1+\beta)  (
M_{\beta}(t_1)\cdots M_{\beta}(t_k)
)$ and (abusing the notation a bit by using the same letter $P$),
%
%
\begin{eqnarray}\label{etoomuchmass}
&& P \bigl( \mu(A) \Gamma(1+\beta) \bigl( M_{\beta}(t_1)
\cdots M_{\beta}(t_k) \bigr)\in B \bigr)
\nonumber\\[-8pt]\\[-8pt]
&&\qquad  > (1+\rho) P \bigl(
(Y_1,\ldots, Y_k) \in B \bigr)\nonumber
\end{eqnarray}
for some $\rho>0$. In fact, since the law of a Mittag--Leffler random
variable is atomless, such a $B$ can be taken to be either a ``SW
corner'' of the type $B=\prod_{j=1}^k (-\infty,x_j]$ for some
$(x_1,\ldots, x_k)\in\mathbb{R}^k$, or its complement.

Choose now $0<\theta<1$ so small that
%
%
\begin{equation}
\label{echoicevep} (1-\theta) (1+\rho)>1
\end{equation}
and
consider the set $A_{k,\theta}$. It follows from (\ref{epresumeconv})
and Hopf's ergodic theorem that
\[
\frac{1}{a_{n_l}} \bigl( S_{\lceil n_l t_1 \rceil} (\mathbf{1} _{A_{k,\theta
}}),\ldots, S_{\lceil n_l t_k \rceil} (\mathbf{1}_{A_{k,\theta}}) \bigr)\Rightarrow
\frac{\mu ( A_{k,\theta} )}{\mu(A)} (Y_1,\ldots, Y_k),
\]
when the law of the
random vector in the left-hand side is still computed with respect to
$\mu_A$. However, since (\ref{eqlonggoal}) holds for
the set $A_{k,\theta}$, we see that
\begin{eqnarray*}
&& P \bigl( (Y_1,\ldots, Y_k) \in B \bigr)
\\
&& \qquad = \lim
_{l\to\infty} \mu_A \biggl( \frac{1}{a_{n_l}} \bigl(
S_{\lceil n_l
t_1 \rceil} (\mathbf{1}_{A_{k,\theta}}), \ldots, S_{\lceil n_l t_k \rceil} (
\mathbf{1}_{A_{k,\theta}}) \bigr) \in \frac{\mu ( A_{k,\theta} )}{\mu(A)} B \biggr)
\\
&&\qquad = \frac{\mu ( A_{k,\theta} )}{\mu(A)} \lim_{l\to\infty} \mu_ {A_{k,\theta}} \biggl(
\frac
{1}{a_{n_l}} \bigl( S_{\lceil n_l
t_1 \rceil} (\mathbf{1}_{A_{k,\theta}}), \ldots,
S_{\lceil n_l t_k \rceil} (\mathbf{1}_{A_{k,\theta}}) \bigr) \in \frac{\mu ( A_{k,\theta} )}{\mu(A)} B
\biggr)
\\
&&\qquad \geq(1-\theta) P \bigl( \mu(A) \Gamma(1+\beta) \bigl( M_{\beta
}(t_1)
\cdots M_{\beta}(t_k) \bigr)\in B \bigr)
\\
&&\qquad > P \bigl( (Y_1,\ldots,Y_k) \in B \bigr),
\end{eqnarray*}
where the last inequality follows from (\ref{etoomuchmass}) and
(\ref
{echoicevep}).
This contradiction shows that, once we prove
(\ref{eupperbound}), this will establish (\ref{eqlonggoal}) for
every measurable $A$ with $0<\mu(A)<\infty$.

We call a nested sequence $(A_0,A_1,\ldots)$ of sets in
(\ref{eapproxDK}) an $\varepsilon$-sequence starting at~$A_0$. Its finite
subsequence $(A_0,A_1,\ldots, A_k)$ will be called an $\varepsilon$-sequence
of length $k+1$ starting at $A_0$ and ending at $A_k$. Let $A$ be a
measurable set with $0<\mu(A)<\infty$. Fix
$0<\theta<1$. Let $0<r<1$ be a small number, to be specified in the
sequel. We construct a nested sequence of sets as follows.

Let $B_0=A$. Construct an $r$-sequence
of length $k+1$ starting at $B_0$, and ending at some set
$B_1\subseteq B_0$. Next, construct an $r^2$-sequence
of length $k+1$ starting at $B_1$,  and ending at some set
$B_2\subseteq B_1$. Proceeding this way we obtain a nested sequence of
measurable sets $A=B_0\supseteq B_1\supseteq B_2\supseteq \cdots,$ such that
\[
\mu(B_n) \geq\prod_{i=1}^n
\bigl(1-r^i\bigr)^k \mu(A),\qquad n=1,2,\ldots.
\]
The sets $(B_n)$ decrease to some set $A_{k,\theta}$ with
\[
\mu(A_{k,\theta}) \geq\prod_{i=1}^\infty
\bigl(1-r^i\bigr)^k \mu(A).
\]
Notice that, by choosing $0<r<1$ small enough, we can ensure that
(\ref{ealmostA}) holds. Note, further, that by construction, for
every $d=1,2,\ldots,$
\[
\mu(A_{k,\theta}) \geq f_d \mu(B_d)\qquad
\mbox{with } f_d = \prod_{i=d+1}^\infty
\bigl(1-r^i\bigr)^k.
\]
Clearly, $f_d\uparrow1$ as $d\to\infty$. Starting with the first line
in (\ref{ediscrF}), we see that
\begin{eqnarray*}
&& F_{n,A_{k,\theta}}(\theta_1 \cdots\theta_k)
\\
&&\qquad \leq
\bigl(1+o(1) \bigr) \frac{1}{a_n^k}
\\
&&\quad\qquad{}\times  \sum_{m_1=0}^\infty
\cdots\sum_{m_k=0}^\infty e^{-n^{-1}\sum_{j=1}^k v_j m_j}
\mu_{B_d} \Biggl( \bigcap_{j=1}^k
T^{- (m_1+\cdots+m_j) } B_d \Biggr) \frac{\mu(B_d)}{\mu ( A_{k,\theta} )}
\\
&&\qquad \leq \bigl(1+o(1) \bigr) \frac{1}{f_d} \frac{1}{a_n^k}
\bigintsss_{B_d} \Biggl[ \Biggl( \sum_{m_1=0}^\infty
\widehat T^{m_1} \mathbf{1}_{B_{d-1}} e^{- v_1m_1/n} \Biggr)
\\
&&\hspace*{127pt}{}\times
\prod_{j=2}^k \Biggl( \sum
_{m_j=0}^\infty\mathbf{1}_{B_d}\circ
T^{m_2+\cdots
+m_j} e^{-
v_jm_j/n} \Biggr) \Biggr] \,d\mu_{B_d}.
\end{eqnarray*}
Using repeatedly uniform convergence as in
(\ref{euseDK}) above, we conclude, as in the case of the
corresponding lower bound calculation that
\begin{eqnarray*}
&& F_{n,A_{k,\theta}}(\theta_1 \cdots\theta_k)
\\
&&\qquad \leq
\bigl(1+o(1) \bigr) \frac{1}{f_d} \frac{1}{a_n^k} \frac{\mu ( B_{d-1}  )v_1}{n}
\\
&&\quad\qquad{} \times\sum
_{i=0}^{\infty} a_i
e^{-
v_1i/n}\bigintsss_{B_d} \Biggl[ \Biggl( \sum
_{m_2=0}^\infty \widehat T^{m_2}
\mathbf{1}_{B_{d-1}} e^{- v_2m_2/n} \Biggr)
\\
&&\hspace*{124pt}{}\times \prod
_{j=3}^k \Biggl( \sum_{m_j=0}^\infty
\mathbf{1}_{B_d}\circ T^{m_3+\cdots
+m_j} e^{-
v_jm_j/n} \Biggr) \Biggr]
\,d\mu_{B_d}
\\
&&\qquad \leq \cdots\leq \bigl( 1 + o(1) \bigr) \frac{1}{f_d} \biggl(
\frac
{\mu
(B_{d-1})}{na_n} \biggr)^k ( v_1 \cdots v_k )
\prod_{j=1}^k \Biggl( \sum
_{i=0}^{\infty} a_i e^{-v_j i /n} \Biggr)
\\
&&\qquad \leq \bigl( 1 + o(1) \bigr) \frac{1}{f_d f_{d-1}^k} \biggl( \frac
{\mu(A_{k,\theta})}{na_n}
\biggr)^k ( v_1 \cdots v_k ) \prod
_{j=1}^k \Biggl( \sum_{i=0}^{\infty}
a_i e^{-v_j i /n} \Biggr).
\end{eqnarray*}
As in the case of the lower bound, Karamata's Tauberian theorem shows
that
\begin{eqnarray*}
F_{n,A_{k,\theta}}(\theta_1 \cdots\theta_k)
&\leq&  \bigl(
1+o(1) \bigr) \frac{1}{f_d f_{d-1}^k} \bigl( \mu(A_{k,\theta}) \Gamma(1+\beta)
\bigr)^k \prod_{j=1}^k
\frac{a(n/v_j)}{a_n}
\\
&\to&\frac{1}{f_d f_{d-1}^k} \bigl( \mu(A_{k,\theta})\Gamma(1+\beta)
\bigr)^k \prod_{j=1}^k
v_j^{-\beta}
\end{eqnarray*}
as $n\to\infty$. Since this is true for every $d\geq1$, we can let
now $d\to\infty$ to obtain~(\ref{ealmostA}), and the proof of the
theorem is complete.
\end{pf}

%
\begin{remark} \label{rkalsoconv}
It follows immediately from Theorem~\ref{tgenDK} and continuity of
the limiting Mittag--Leffler process that for the continuous process
$(\widetilde S_n)$
defined in~(\ref{einterpprocess}), strong distributional convergence
as in (\ref{eqgoalLemma1}) also holds, either in $D[0,\infty)$ or
in $C[0,\infty)$.
\end{remark}

We use the strong distributional convergence obtained in Theorem
\ref{tgenDK} in the following proposition.

%
\begin{proposition} \label{plineartime}
Under the assumptions of Theorem
\ref{tgenDK}, let $A$ be a Darling--Kac set with $0<\mu(A)<\infty$,
such that (\ref{eqstrongDK}) is satisfied, and suppose that the
function $f$ is supported by $A$. Define a
probability measure on $E$ by $\mu_n (\cdot) = \mu(\cdot\cap
\{\varphi\leq n\}) / \mu(\{\varphi\leq n \})$, where $\varphi$ is
the first entrance time of $A$. Let $0 \leq t_1 < \cdots< t_H$, $H
\geq1$, and\vspace*{1pt} fix $L \in\mathbb{N}$ with $t_H \leq L$. Then under
$\mu_{nL}$, the sequence $ ( S_{\lceil nt_h \rceil}(f) /
a_n )_{h=1}^H$ converges\vspace*{1pt}
weakly in $\mathbb{R}^H$ to the random vector $ (\mu(f) \Gamma
(1+\beta)
M_{\beta}(t_h-T_{\infty}^{(L)})_+ )_{h=1}^H$, where
$T_{\infty}^{(L)}$ is a random variable independent of the
Mittag--Leffler process $M_\beta$, with $P ( T_{\infty
}^{(L)}\leq
x )=(x/L)^{1-\beta}$, $ 0\leq x\leq L$.
\end{proposition}

\begin{pf}
Since $T$ preserves measure $\mu$, for the duration of the proof we
may and will modify the definition of $S_n$ to $S_n(f) =
\sum_{k=0}^{n-1} f \circ T^k$, $n=1,2,\ldots.$  Fix $\theta_1, \ldots,
\theta_H \in\mathbb{R}$ and let $\lambda\in\mathbb{R}$. Since $f$
is supported by $A$, we have, as $n \to\infty$,
\begin{eqnarray*}
&& \mu_{n L} \Biggl( \frac{1}{a_n} \sum
_{h=1}^H \theta_h S_{\lceil n t_h
\rceil}
(f) > \lambda \Biggr)
\\
&&\qquad \sim\mu_{n L} \Biggl( A^c \cap \Biggl\{
\frac{1}{a_n} \sum_{h=1}^H
\theta_h S_{\lceil n t_h \rceil} (f) > \lambda \Biggr\} \Biggr)
\\
&&\qquad = \mu(\varphi\leq nL)^{-1} \sum_{m=1}^{nL}
\mu \Biggl( A_m \cap \Biggl\{ \frac{1}{a_n} \sum
_{h=1}^H \theta_h S_{\lceil n t_h \rceil}
(f) > \lambda \Biggr\} \Biggr)
\\
&&\qquad \sim\mu(\varphi\leq nL)^{-1} \sum_{m=1}^{nL}
\mu \Biggl( A_m \cap T^{-m} \Biggl\{ \frac{1}{a_n}
\sum_{h=1}^H \theta_h
S_{(\lceil n t_h
\rceil- m)_+} (f) > \lambda \Biggr\} \Biggr)
\\
&&\qquad = \int_A \frac{1}{\mu(\varphi\leq nL)} \sum
_{m=1}^{nL} \widehat{T}{}^m
\mathbf{1}_{A_m} \cdot\mathbf{1}_{\{ \sum_{h=1}^H \theta_h
S_{(\lceil n t_h
\rceil
- m)_+} (f) > \lambda a_n \}} \,d\mu.
\end{eqnarray*}

Note that the measure on $E$ defined by $\eta(\cdot) = \int_{\cdot} K
\,d\mu$ with $K$ in (\ref{eqstrongDK}) is necessarily a probability
measure. We conclude by (\ref{eqstrongDK}) that
%
%
\begin{eqnarray}\label{eqprodexpression}
&& \mu_{n L} \Biggl( \frac{1}{a_n} \sum
_{h=1}^H \theta_h S_{\lceil n t_h\rceil}
(f) > \lambda \Biggr)
\nonumber\\[-8pt]\\[-8pt]
&&\qquad  \sim\sum_{m=1}^{nL}
\eta \Biggl( \frac{1}{a_n} \sum_{h=1}^H
\theta_h S_{(\lceil nt_h \rceil- m)_+} (f) > \lambda \Biggr) p_n(m), \nonumber
\end{eqnarray}
where $p_n(j) = \mu(A_j) / \sum_{m=1}^{nL} \mu(A_m)$, $j=1, \ldots,
nL$, is a probability mass function. Let $T_n^{(L)}$ be a discrete
random variable with this probability mass
function, independent of $S_{\lceil n \cdot\rceil}(f)$, which is, in
turn, governed by the probability measure $\eta$. If we declare that
$T_n^{(L)}$ is defined on some probability space $ ( \Omega_n,
{\mathcal F}_n, P_n )$, then the right-hand side of
(\ref{eqprodexpression}) becomes
\[
( \eta\times P_n ) \Biggl( \frac{1}{a_n} \sum
_{h=1}^H \theta_h S_{(\lceil nt_h \rceil-T_n^{(L)} )_+}
(f) > \lambda \Biggr).
\]

Since $\eta$ is a probability measure absolutely continuous with
respect to $\mu$, it follows from the strong distributional
convergence in Theorem~\ref{tgenDK} that
%
%
\begin{equation}
\frac{1}{a_n} S_{\lceil n \cdot\rceil}(f) \Rightarrow\mu(f)\Gamma (1+\beta)
M_{\beta}(\cdot) \qquad\mbox{in } D[0,L], \label{eqfromLemma1}
\end{equation}
when the law in the left-hand side is computed with respect to
$\eta$. On the other hand, by the regular variation of the wandering
rate sequence and (\ref{eqwanderingmu}), for $x \in[0,L]$,
%
%
\begin{eqnarray}
P_n \biggl( \frac{T_n^{(L)}}{n} \leq x \biggr) = \sum
_{m=1}^{\lceil nx
\rceil} p_n(m) \sim\frac{w_{\lceil nx \rceil}}{w_{nL}}
\sim \biggl( \frac{x}{L} \biggr)^{1-\beta}, \label{eqbeta}
\end{eqnarray}
which is precisely the law of $T_{\infty}^{(L)}$.
We can put together (\ref{eqfromLemma1}), (\ref{eqbeta}), and
independence between $S_n $ and $T_n^{(L)} $ to obtain
\begin{eqnarray*}
&& \mu_{nL} \Biggl( \frac{1}{a_n} \sum_{h=1}^H
\theta_h S_{\lceil nt_h
\rceil}(f) > \lambda \Biggr) 
\\
&&\qquad \to P \Biggl(
\mu(f)\Gamma(1+\beta) \sum_{h=1}^H
\theta_h M_{\beta}\bigl(\bigl(t_h -
T_{\infty}^{(L)}\bigr)_+\bigr) > \lambda \Biggr)
\end{eqnarray*}
for all continuity points $\lambda$ of the right-hand side, and all
$\theta_1 \cdots\theta_H \in\mathbb{R}$ by, for example, Theorem
13.2.2 in
\citet{whitt2002}. This proves the proposition.
\end{pf}

\section{Proof of the main theorem} \label{secproof}

In this section, we prove Theorem~\ref{tFCLT}. We start with several
preliminary results. The first lemma explains the asymptotic relation
(\ref{eqasycn}).

%
\begin{lemma} \label{lmomentgrowth}
Under the assumptions of Proposition~\ref{plineartime}, assume,
additionally, that the set $A$ supporting $f$ is a Darling--Kac set. Let
$0<\alpha<2$. If $1<\alpha<2$, assume, additionally, that $f\in
L^2(\mu)$, and that either:
\begin{longlist}[(ii)]
\item[(i)] $A$ is a uniform set for $|f|$, or

\item[(ii)] $f$ is bounded.
\end{longlist}

Then
%
%
\begin{equation}\label{eqasynormalizing}
\biggl( \int_E \bigl|S_n(f)\bigr|^{\alpha} \,d\mu
\biggr)^{1/\alpha} \sim \bigl|\mu(f)\bigr| C_{\alpha,\beta} a_n
w_n^{1/\alpha} \qquad\mbox{as } n \to \infty
\end{equation}
and (\ref{eqasycn}) holds.
\end{lemma}

\begin{pf}
It is an elementary calculation to check that
(\ref{eqasynormalizing}) implies (\ref{eqasycn}), so in the sequel
we concentrate on checking (\ref{eqasynormalizing}). It follows from
(\ref{eqwanderingmu}) and the fact that $f$ is supported by $A$
that
%
%
\begin{equation}
\biggl( \int_E \bigl|S_n(f)\bigr|^{\alpha} \,d\mu
\biggr)^{1/\alpha} = a_n \bigl(\mu(\varphi\leq n)
\bigr)^{1/\alpha} A_n^{(\alpha)} \sim a_n
w_n^{1/\alpha} A_n^{(\alpha)}, \label{eqexpressioncn}
\end{equation}
where $A_n^{(\alpha)} =( \int_E |S_n(f) / a_n |^{\alpha} \,d \mu_n
)^{1/\alpha}$. Therefore, proving (\ref{eqasynormalizing}) reduces
to checking that
%
%
\begin{equation}
\label{emomentconv} A_n^{(\alpha)} \to\bigl|\mu(f)\bigr| C_{\alpha,\beta}
\qquad\mbox{as }n\to\infty.
\end{equation}
If $\alpha=1$ and $f$ is nonnegative, then this
follows by direct calculation, using the definition of
$C_{\alpha,\beta}$. If $f$ is not necessarily nonnegative, we can use
the obvious bound $-S_n(|f|)\leq S_n(f)\leq S_n(|f|)$ together with
the so-called Pratt lemma; see \citet{pratt1960}, or Problem 16.4(a)
in \citet{billingsley1995}.

It remains to consider the case $\alpha\in(0,1)\cup(1,2)$.
Proposition~\ref{plineartime} shows that $ (
A_n^{(\alpha)} )$ is the sequence of the $\alpha$-norms of a
weakly converging sequence, and the expression in the right-hand side
of (\ref{emomentconv}) is easily seen to be the $\alpha$-norm of
the weak
limit. Therefore, our statement will follow once we show that this
weakly convergent sequence is uniformly integrable, which we proceed
now to do.

Suppose first that $0<\alpha<1$. Recalling the relation
(\ref{eqprop387}) and the fact that $T$ preserves measure $\mu$, we
see that
%
%
\begin{eqnarray}\label{eunifl1}
\sup_{n \geq1} \int_E \biggl\llvert
\frac{S_n(f)}{a_n} \biggr\rrvert \,d\mu_n &=& \sup_{n \geq1}
\frac{1}{a_n \mu(\varphi\leq n)} \int_E \bigl|S_n(f)\bigr| \,d\mu
\nonumber\\[-8pt]\\[-8pt]
&\leq&\sup_{n \geq1} \frac{n}{a_n \mu(\varphi\leq n)} \int_E
|f| \,d\mu <\infty,\nonumber
\end{eqnarray}
which proves uniform integrability in this case.

Finally, we consider the case $1<\alpha<2$, when it is sufficient to
prove that
%
%
\begin{equation}
\sup_{n \geq1} \int_E \biggl(
\frac{S_n(f)}{a_n} \biggr)^2 \,d\mu_n < \infty.
\label{equnifintegral}
\end{equation}
Under the assumption (i), since $f$ is supported by $A$, we can use
the duality relation (\ref{edualrel}) to write
\begin{eqnarray*}
\int_E S_n(f)^2 \,d\mu&=& n \int
_E f^2 \,d\mu+ \sum
_{k=1}^n \sum_{l=1, k
\neq l}^n
\int_E f \circ T^k f \circ T^l d
\mu
\\
&=& n \int_E f^2 \,d\mu+ 2\sum
_{k=1}^{n-1} \sum_{j=1}^{n-k}
\int_A \widehat{T}{}^j f \cdot f \,d\mu,
\end{eqnarray*}
so that
\begin{eqnarray*}
&& \int_E \biggl( \frac{S_n(f)}{a_n} \biggr)^2
\,d\mu_n
\\
&&\qquad \leq\frac{n}{a_n^2
\mu(\varphi\leq n)} \int_E
f^2 \,d\mu+ \frac{2}{a_n^2 \mu(\varphi
\leq n)} \sum_{k=1}^{n-1}
\sum_{j=1}^{n-k} \int_A
\widehat{T}{}^j |f| \cdot|f| \,d\mu.
\end{eqnarray*}
Clearly, $n/ ( a_n^2 \mu(\varphi\leq n) ) \to0$. Further,
since $A$ is uniform for $|f|$,
\begin{eqnarray*}
&& \frac{1}{a_n^2 \mu(\varphi\leq n)} \sum_{k=1}^{n-1} \sum
_{j=1}^{n-k} \int_A
\widehat{T}{}^j |f| \cdot|f| \,d\mu
\\
&&\qquad \leq\frac{n}{a_n \mu(\varphi\leq n)} \int
_A \frac{1}{a_n} \sum_{j=1}^n
\widehat{T}{}^j |f| \cdot|f| \,d\mu
\sim\mu\bigl(|f|\bigr)^2 \frac{n}{a_n \mu(\varphi\leq n)}.
\end{eqnarray*}
Using (\ref{eqprop387}), we see that (\ref{equnifintegral})
follows. On the other hand, under the assumption~(ii), the ratio
$S_n(f) / S_n(\mathbf{1}_A)$ is bounded, hence for some finite
$C> 0$,
\[
\sup_{n \geq1} \int_E \biggl(
\frac{S_n(f)}{a_n} \biggr)^2 \,d\mu_n \leq C\sup
_{n \geq1} \int_E \biggl(
\frac{S_n(\mathbf{1}_A)}{a_n} \biggr)^2 \,d\mu_n.
\]
However, the Darling--Kac property of $A$ means that it is uniform for
$\mathbf{1}_A$, and so we are, once again, under the assumption (i).
\end{pf}

In preparation for the proof of Theorem~\ref{tFCLT}, we introduce a
useful decomposition of the process $\mathbf{X}$ given in
(\ref{etheprocess}). We begin by decomposing the local L\'evy
measure $\rho$ into a sum of two parts, corresponding to ``large
jumps'' and ``small jumps.'' Let
\begin{eqnarray*}
\rho_1 (\cdot) &=& \rho \bigl(\cdot\cap\bigl\{|x| > 1 \bigr\} \bigr),
\\
\rho_2 (\cdot) &=& \rho \bigl(\cdot\cap\bigl\{|x| \leq1 \bigr\} \bigr)
\end{eqnarray*}
and let $M_1,  M_2$ be independent homogeneous symmetric infinitely
divisible random
measures, without a Gaussian component, with the same control measure
$\mu$ and local L\'evy measures $\rho_1,  \rho_2$ accordingly. Under
the integrability assumptions (\ref{eqintegrabilitycond}), the
stochastic processes $X_n^{(i)} = \int_E f \circ T^n(x) \,dM_i(x),
n=1,2,\ldots,$ for $i=1,2$, are independent stationary infinitely
divisible processes,
and $X_n = X_n^{(1)} + X_n^{(2)},   n=1,2,\ldots.$

Our final lemma shows that, from the point of view of the central limit
behavior in the case $0<\alpha<1$, the contribution of the process
$ (X_n^{(2)} )$, corresponding to the ``small jumps,'' is negligible.

%
\begin{lemma} \label{lsmalljumps}
If $0<\alpha<1$, then
%
%
\begin{equation}
\frac{1}{c_n} \sum_{k=1}^n
X_k^{(2)} \stackrel{p} {\to} 0. \label{eqnegp}
\end{equation}
\end{lemma}

\begin{pf} By Chebyshev's inequality, for any $\epsilon>0$,
\[
P \Biggl( \Biggl|\sum_{k=1}^n
X_k^{(2)}\Biggr| > \epsilon c_n \Biggr) \leq
\frac{n}{\epsilon c_n} E\bigl|X_1^{(2)}\bigr| \to0
\]
(since $c_n \in RV_{\beta+ (1-\beta)/\alpha}$ implies $n/c_n \to0$ in
the case $0<\alpha<1$) as long as the expectation $E|X_1^{(2)}|$ is
finite. Since for every $p_1>p_0$ in (\ref{eqorginregularity}) and
$p_1\geq1$,
\[
\int_E\int_\mathbb{R}\bigl|xf(s)\bigr|\mathbf{1}
\bigl( \bigl|xf(s)\bigr|>1 \bigr) \rho_2(dx) \mu(ds) \leq\int
_{-1}^1 |x|^{p_1} \rho(dx) \int
_E \bigl|f(s)\bigr|^{p_1} \mu (ds),
\]
the expectation is finite because, by (\ref{eqintegrabilitycond}), we
can find $p_1$ as above such that $\int_E |f|^{p_1}  \,d\mu<\infty$.
\end{pf}

\begin{pf*}{Proof of Theorem~\ref{tFCLT}}
We start with proving the finite-dimensional weak convergence, for
which it enough to show the convergence
\[
\frac{1}{c_n} \sum_{h=1}^H
\theta_h \sum_{k=1}^{\lceil nt_h \rceil}
X_k \Rightarrow\bigl|\mu(f)\bigr|\sum_{h=1}^H
\theta_h Y_{\alpha,
\beta}(t_h)
\]
for all $H \geq1$, $0 \leq t_1 < \cdots< t_H$, and $\theta_1 \cdots
\theta_H \in\mathbb{R}$. Conditions for weak convergence of
infinitely divisible random
variables [see, e.g., Theorem 15.14 in \citet{kallenberg2002}]
simplify in this one-dimensional symmetric case to
%
%
\begin{eqnarray} \label{eqfirstcond}
&& \int_E \Biggl( \frac{1}{c_n} \sum
_{h=1}^H \theta_h
S_{\lceil nt_h
\rceil}(f) \Biggr)^2 \int
_0^{rc_n / |\sum\theta_h
S_{\lceil nt_h
\rceil}(f)|}
x \rho(x,\infty) \,dx \,d\mu\nonumber
\\
&&\qquad \to\frac{r^{2-\alpha} C_\alpha}{2-\alpha} \bigl|\mu(f)\bigr|^\alpha
\\
&&\quad\qquad{} \times \int_{[0,\infty)}\int_{\Omega^{\prime}} \Biggl\llvert \sum_{h=1}^H
\theta_h M_{\beta}\bigl((t_h-x)_+,
\omega^{\prime}\bigr) \Biggr\rrvert ^{\alpha} P^{\prime}\bigl(d
\omega^{\prime}\bigr) \nu(dx)\nonumber
\end{eqnarray}
and
%
%
\begin{eqnarray}\label{eqomit}
&& \int_{E} \rho \Biggl( rc_n \Biggl| \sum
_{h=1}^H \theta_h S_{\lceil nt_h
\rceil}
(f)\Biggr|^{-1}, \infty \Biggr) \,d\mu\nonumber
\\
&&\qquad \to r^{-\alpha} C_\alpha \bigl|\mu(f)\bigr|^\alpha
\\
&&\quad\qquad {}\times \int
_{[0,\infty)} \int_{\Omega^{\prime}} \Biggl\llvert \sum
_{h=1}^H \theta_h
M_{\beta}\bigl((t_h-x)_+,\omega^{\prime}\bigr) \Biggr
\rrvert ^{\alpha} P^{\prime}\bigl(d\omega^{\prime}\bigr) \nu(dx)\nonumber
\end{eqnarray}
for every $r>0$. Fix $L
\in\mathbb{N}$ with $t_H \leq L$ and $r>0$.

Since the argument for (\ref{eqfirstcond}) and
the argument for (\ref{eqomit}) are very similar, we only
prove (\ref{eqfirstcond}). By Proposition
\ref{plineartime} and Skorohod's embedding theorem, there is some
probability space $ ( \Omega^*, {\mathcal F}^*, P^* )$ and
random variables $Y$, $Y_n$, $n=1,2,\ldots$ defined on that space such
that, for every $n$, the law of $Y_n$ coincides with the law of
$a_n^{-1}\sum_{h=1}^H \theta_h S_{\lceil nt_h \rceil}(f)$ under
$\mu_{nL}$, the law of $Y$ coincides with the law of $\mu(f)\Gamma
(1+\beta)
\sum_{h=1}^H \theta_h M_{\beta}((t_h - T_{\infty}^{(L)})_+)$ under
$P^{\prime}$, and $Y_n\to Y$ $P^*$-a.s.

Introduce a function
\[
\psi(y) = y^{-2} \int_0^{ry} x
\rho(x,\infty) \,dx,\qquad y>0,
\]
so that the expression in the left-hand side of (\ref{eqfirstcond})
becomes
\[
\int_E \psi \biggl( \frac{c_n}{|\sum_{h=1}^H \theta_h S_{\lceil
nt_h \rceil}(f)|} \biggr) \,d\mu =
\mu (\varphi\leq nL ) E^* \biggl[ \psi \biggl( \frac{c_n}{a_n |Y_n|} \biggr) \biggr].
\]
By Karamata's theorem [see, e.g., Theorem 0.6 in \citet{resnick1987}],
\[
\psi(y) \sim\frac{r^2}{2-\alpha} \rho(ry,\infty) \qquad\mbox{as } y \to\infty,
\]
so that, as $n\to\infty$,
%
%
\begin{eqnarray}\label{eqnegligbleuniform}
\qquad&& \mu (\varphi\leq nL ) \psi \biggl( \frac{c_n}{a_n|Y_n|} \biggr) \nonumber
\\
&&\qquad \sim\frac{r^2}{2-\alpha} \mu (\varphi\leq nL ) |Y_n|^{\alpha}
\rho \bigl( rc_na_n^{-1},\infty \bigr)
\\
&&\quad\qquad{} + \frac{r^{2}}{2-\alpha} \mu (\varphi\leq nL ) \rho \bigl(rc_n
a_n^{-1}, \infty \bigr) \biggl( \frac{\rho (rc_n a_n^{-1}
|Y_n|^{-1},\infty )}{\rho (rc_n a_n^{-1},\infty )} -
|Y_n|^{\alpha} \biggr).\nonumber
\end{eqnarray}
By (\ref{eqasycn}), Lemma~\ref{lmomentgrowth} and
(\ref{eqwanderingmu}),
%
%
\begin{equation}
\qquad\rho\bigl(rc_n a_n^{-1},\infty\bigr) \sim
r^{-\alpha} C_\alpha \bigl( \Gamma(1+\beta) \bigr)^{-\alpha}
\bigl(\mu(\varphi\leq n) \bigr)^{-1} \qquad\mbox{as } n \to\infty.
\label{eqrhoandmu}
\end{equation}
This, together with the basic properties of regularly varying
functions of a negative index [see, e.g., Proposition 0.5,
\citet{resnick1987}], shows that the second term in the right-hand
side of (\ref{eqnegligbleuniform}) converges to $0$. Therefore,
\[
\mu (\varphi\leq nL ) \psi \biggl( \frac{c_n}{a_n|Y_n|} \biggr) \to
\frac{r^{2-\alpha}}{2-\alpha} C_\alpha L^{1-\beta
} \biggl( \frac{|Y|}{\Gamma(1+\beta)}
\biggr)^{\alpha}.
\]
Integrating the limit yields
\begin{eqnarray*}
&& E^* \biggl[ \frac{r^{2-\alpha}}{2-\alpha} C_\alpha L^{1-\beta
} \biggl(
\frac{|Y|}{\Gamma(1+\beta)} \biggr)^{\alpha} \biggr]
\\
&&\qquad = \frac{r^{2-\alpha}}{2-\alpha}
C_\alpha L^{1-\beta} \bigl|\mu (f)\bigr|^\alpha E^{\prime}
\Biggl[ \sum_{h=1}^H \theta_h
M_{\beta}\bigl(\bigl(t_h-T_{\infty}^{(L)}
\bigr)_+\bigr) \Biggr]^{\alpha}
\\
&&\qquad = \frac{r^{2-\alpha}C_\alpha}{2-\alpha} \bigl|\mu(f)\bigr|^\alpha \int_{[0,\infty)}
\int_{\Omega^{\prime}} \Biggl( \sum_{h=1}^H
\theta_h M_{\beta}\bigl((t_h-x)_+,
\omega^{\prime
}\bigr) \Biggr)^{\alpha} P^{\prime}\bigl(d
\omega^{\prime}\bigr) \nu(dx),
\end{eqnarray*}
which is exactly the right-hand side of
(\ref{eqfirstcond}). Therefore, in order to complete the proof of
(\ref{eqfirstcond}), we only need to justify taking the limit inside
the integral. For this purpose, we use, once again, Pratt's lemma. We
need to exhibit random variables~$G_n$, $n=0,1,2,\ldots$ on $ (
\Omega^*, {\mathcal F}^*, P^* )$ such that
%
%
\begin{eqnarray}
\mu (\varphi\leq nL ) \psi \biggl( \frac{c_n}{a_n|Y_n|} \biggr) &\leq&
G_n, \qquad P^*\mbox{-a.s.}, \label{eq1stPratt}
\\
G_n &\to& G_0, \qquad P^*\mbox{-a.s.}, \label{eqeasyPratt}
\\
E^* G_n &\to& E^* G_0 \in[0,\infty). \label{eq3rdPratt}
\end{eqnarray}

We start with writing [using (\ref{eqrhoandmu})]
\begin{eqnarray*}
&& \mu (\varphi\leq nL ) \psi \biggl( \frac{c_n}{a_n|Y_n|} \biggr)
\\
&&\qquad  \leq
C_1 \frac{\psi ( c_n a_n^{-1} |Y_n|^{-1}  )}{\psi(c_n
a_n^{-1})} \mathbf{1}_{\{ c_n > a_n |Y_n| \}} + C_1
\frac{\psi ( c_n a_n^{-1} |Y_n|^{-1}  )}{\psi(c_n
a_n^{-1})} \mathbf{1}_{\{ c_n \leq a_n |Y_n| \}},
\end{eqnarray*}
where $C_1>0$ is a constant.
Suppose first that $1 \leq\alpha< 2$, and choose
$0<\xi<2-\alpha$. Then by the Potter bounds [see
Proposition 0.8 in \citet{resnick1987}], for some constant $C_2 > 0$,
\[
\frac{\psi ( c_n a_n^{-1} |Y_n|^{-1}  )}{\psi(c_n
a_n^{-1})} \mathbf{1}_{\{ c_n > a_n |Y_n| \}} \leq C_2
\bigl(|Y_n|^{\alpha-
\xi} + |Y_n|^{\alpha+ \xi}\bigr)
\]
for all $n$ large enough.
Further, since $y^2 \psi(y) \to0$ as $y \downarrow0$, we have, for
some constant $C_3>0$,
\[
\frac{\psi ( c_n a_n^{-1} |Y_n|^{-1}  )}{\psi(c_n
a_n^{-1})} \mathbf{1}_{\{ c_n \leq a_n |Y_n| \}} \leq C_3 \biggl(
\frac{a_n}{c_n} \biggr)^2 \frac{|Y_n|^2}{\psi(c_n
a_n^{-1})},
\]
hence, for
some constant $C_4>0$,
%
%
\begin{equation}
\qquad\quad \mu (\varphi\leq nL ) \psi \biggl( \frac{c_n}{a_n|Y_n|} \biggr) \leq
C_4 \biggl( |Y_n|^{\alpha- \xi} +
|Y_n|^{\alpha+ \xi} + \biggl( \frac{a_n}{c_n}
\biggr)^2 \frac{|Y_n|^2}{\psi(c_n a_n^{-1})} \biggr) \label{eqmorethan1}
\end{equation}
for all $n$ (large enough) and all realizations. We take
\begin{eqnarray*}
G_n &=& C_4 \biggl( |Y_n|^{\alpha- \xi} +
|Y_n|^{\alpha+ \xi} + \biggl( \frac{a_n}{c_n}
\biggr)^2 \frac{|Y_n|^2}{\psi(c_n a_n^{-1})} \biggr), \qquad n=1,2,\ldots,
\\
G_0 &=& C_4\bigl(|Y|^{\alpha- \xi} +
|Y|^{\alpha+ \xi}\bigr).
\end{eqnarray*}
Then (\ref{eq1stPratt}) holds by construction, while
(\ref{eqeasyPratt}) follows from the fact that
\[
\biggl( \frac{a_n}{c_n} \biggr)^2 \frac{1}{\psi(c_n a_n^{-1})} \in
RV_{(1-\beta)(1-2/\alpha)}
\]
and $(1-\beta)(1-2/\alpha) < 0$. Keeping this in mind, and recalling
that, by (\ref{equnifintegral}) (which holds also for $\alpha=1$
under the assumptions of the theorem), $\sup_{n\geq1}E^*Y_n^2<\infty$,
we
obtain the uniform integrability implying (\ref{eq3rdPratt}). This
proves (\ref{eqfirstcond}) in the case $1\leq\alpha<2$.

If $0 < \alpha< 1$, then Lemma~\ref{lsmalljumps} allows us to
assume, without loss of generality, that $\rho(x\dvtx |x| \leq1) = 0$.
Then $\psi$ is bounded on $(0,1]$, so that for some $C_5>0$,
\[
\frac{\psi(c_n a_n^{-1} |Y_n|^{-1})}{\psi(c_n a_n^{-1})} \mathbf {1}_{\{ c_n
\leq a_n |Y_n| \}} \leq C_5
\frac{a_n}{c_n} \frac{|Y_n|}{\psi(c_n
a_n^{-1})}
\]
and the upper bound (\ref{eqmorethan1}) is replaced with
\[
\mu (\varphi\leq nL ) \psi \biggl( \frac{c_n}{a_n
|Y_n|} \biggr) \leq
C_6 \biggl( |Y_n|^{\alpha-
\xi} + |Y_n|^{\alpha+ \xi}
+ \frac{a_n}{c_n} \frac{|Y_n|}{\psi(c_n a_n^{-1})} \biggr)
\]
for some $C_6>0$, where we now choose $0<\xi<1-\alpha$. Since
\[
\frac{a_n}{c_n} \frac{1}{\psi(c_n a_n^{-1})} \in RV_{(1-\beta
)(1-1/\alpha)}
\]
with $(1-\beta)(1-1/\alpha) < 0$ and $\sup_{n\geq1}E^*|Y_n|<\infty$
by (\ref{eunifl1}), an argument similar to the case $1\leq\alpha<2$
applies here as well. A similar argument proves, in the case
$0<\alpha<1$, the ``positive'' version described in Remark
\ref{rkpositive}.

It remains to prove that the laws in the left-hand side of
(\ref{emain}) are tight in $D[0,L]$ for any fixed $L>0$.
By Theorem 13.5 of \citet{billingsley1999}, it is
enough to show that there exist $\gamma_1 > 1$, $\gamma_2 \geq0$ and
$B>0$ such that
\[
P \Biggl[ \min \Biggl( \Biggl|\sum_{k=1}^{\lceil ns \rceil}
X_k - \sum_{k=1}^{\lceil nr
\rceil}
X_k \Biggr|, \Biggl|\sum_{k=1}^{\lceil nt \rceil}
X_k - \sum_{k=1}^{\lceil ns
\rceil}
X_k \Biggr| \Biggr)\geq\lambda c_n \Biggr] \leq
\frac{B}{\lambda^{\gamma_2}} (t-r)^{\gamma_1}
\]
for all $0 \leq r \leq s \leq t \leq L$, $n \geq1$ and $\lambda>
0$. We start with a simple observation that, in the case
$0<\alpha<1$, we may assume that the function $f$ is bounded. To see
that, note that we can always write $f=f\mathbf{1}_{|f|>M} + f\mathbf
{1}_{|f|\leq
M}$, and use the finite-dimensional convergence in
(\ref{emainpos}) and the fact that $\mu ( f\mathbf
{1}_{|f|>M}
)\to
0$ as $M\to\infty$.

Next, for any $0<\alpha<2$, if $0<t-r < 1/n$, then the probability in
the left-hand side vanishes. If $X_n =
X_n^{(1)} + X_n^{(2)},   n=1,2,\ldots$ be the decomposition described
prior to Lemma~\ref{lsmalljumps}. We start with the part
corresponding to the ``small jumps.'' Note that, by Lemma
\ref{lsmalljumps}, this part is negligible if $0<\alpha<1$ (since we
can apply the lemma to the supremum of the process). Therefore, we
only consider the case $1\leq\alpha<2$,
and prove that there exist $\gamma_1 > 1$,
$\gamma_2 \geq0$ and $B>0$ such that for all $0 \leq s \leq t \leq
L$, $n \geq1$, $|t-s|\geq1/n$ and $\lambda>0$,
%
%
\begin{equation}
\label{epart2} P \Biggl( \Biggl|\sum_{k=1}^{\lceil nt \rceil}
X_k^{(2)} - \sum_{k=1}^{\lceil ns
\rceil}
X_k^{(2)} \Biggr| \geq\lambda c_n \Biggr) \leq
\frac{B}{\lambda^{\gamma_2}} (t-s)^{\gamma_1}.
\end{equation}
Note that L\'evy--It{\^o} decomposition yields
\begin{eqnarray*}
&& \sum_{k=1}^{\lceil nt \rceil}X_k^{(2)}
- \sum_{k=1}^{\lceil ns
\rceil}X_k^{(2)}
\\
&&\qquad  \stackrel{d} {=} \int_E S_{\lceil nt \rceil-
\lceil ns \rceil} (f)
\,dM_2
\\
&&\qquad \stackrel{d} {=} \iint_{|x S_{\lceil nt \rceil- \lceil ns \rceil
} (f)|
\leq\lambda c_n} x S_{\lceil nt \rceil- \lceil ns
\rceil} (f)
\,d\overline{N}_2
\\
&&\quad\qquad{} + \iint_{|x S_{\lceil nt \rceil-
\lceil ns
\rceil} (f)| > \lambda c_n} x S_{\lceil nt \rceil-
\lceil ns \rceil}
(f) \,dN_2,
\end{eqnarray*}
where $N_2$ is a Poisson random measure on $\mathbb{R} \times E$ with
mean measure $\rho_2 \times\mu$ and $\overline{N}_2 \equiv N_2 -  (
\rho_2
\times\mu )$. Therefore,
%
%
\begin{eqnarray}\label{eqLevyIto}
&&P \Biggl( \Biggl|\sum_{k=1}^{\lceil nt \rceil}X_k^{(2)}
- \sum_{k=1}^{\lceil ns \rceil}X_k^{(2)}
\Biggr| \geq\lambda c_n \Biggr)
\nonumber
\\
&&\qquad \leq P \biggl( \biggl| \iint_{|x S_{\lceil nt \rceil- \lceil ns
\rceil}
(f)| \leq\lambda c_n}   x S_{\lceil nt \rceil-
\lceil ns \rceil} (f) \,d
\overline{N}_2 \biggr| \geq\lambda c_n \biggr)
\\
&&\quad\qquad{} + P \biggl( \biggl|
\iint_{|x S_{\lceil nt \rceil- \lceil ns \rceil}
(f)| >
\lambda c_n}  x S_{\lceil nt \rceil- \lceil ns
\rceil} (f) \,dN_2 \biggr| >
0 \biggr). \nonumber
\end{eqnarray}
It follows from (\ref{eqorginregularity}) that for some constant
$C_1>0$,
\begin{eqnarray*}
&& P \biggl( \biggl| \iint_{|x S_{\lceil nt \rceil- \lceil ns \rceil}
(f)| \leq\lambda c_n}  x S_{\lceil nt \rceil-
\lceil ns \rceil} (f) \,d
\overline{N}_2 \biggr| \geq\lambda c_n \biggr)
\\
&&\qquad \leq
\frac{1}{\lambda^2 c_n^2} E \biggl\llvert \iint_{|x S_{\lceil
nt \rceil- \lceil ns \rceil} (f)| \leq\lambda c_n} x
S_{\lceil nt \rceil- \lceil ns \rceil} (f) \,d\overline{N}_2 \biggr\rrvert ^2
\\
&&\qquad = \frac{1}{\lambda^2 c_n^2} \iint_{|x S_{\lceil nt \rceil-
\lceil ns \rceil} (f)| \leq\lambda c_n} \bigl\llvert x
S_{\lceil nt \rceil- \lceil ns \rceil} (f) \bigr\rrvert ^2 \rho_2(dx) \,d\mu
\\
&&\qquad \leq4 \int_E \biggl( \frac{S_{\lceil nt \rceil- \lceil ns \rceil}
(f)}{\lambda c_n}
\biggr)^2 \int
_0^{\lambda c_n / |S_{\lceil nt
\rceil- \lceil ns \rceil} (f)|} x
\rho_2(x,\infty) \,dx \,d\mu
\\
&&\qquad \leq\frac{C_1}{ \lambda^{p_0}} \frac{1}{c_n^{p_0}} \int_E
\bigl|S_{\lceil
nt \rceil- \lceil ns \rceil} (f)\bigr|^{p_0} \,d\mu.
\end{eqnarray*}
Similarly, for some constant $C_2>0$,
\begin{eqnarray*}
&& P \biggl( \biggl| \iint_{|x S_{\lceil nt \rceil- \lceil ns \rceil}
(f)| >
\lambda c_n} x S_{\lceil nt \rceil- \lceil ns
\rceil} (f)
\,dN_2\biggr| > 0 \biggr)
\\
&&\qquad \leq P \bigl(N_2 \bigl\{ \bigl|x
S_{\lceil nt \rceil- \lceil ns \rceil} (f)\bigr| > \lambda c_n \bigr\} \geq1 \bigr)
\\
&&\qquad \leq EN_2 \bigl\{ \bigl|x S_{\lceil nt \rceil- \lceil ns \rceil} (f)\bigr| > \lambda
c_n \bigr\}
\\
&&\qquad = 2 \int_E \rho_2 \bigl(\lambda
c_n \bigl|S_{\lceil nt \rceil- \lceil ns
\rceil} (f)\bigr|^{-1},\infty \bigr) \,d\mu
\\
&&\qquad \leq\frac{C_2}{ \lambda^{p_0}} \frac{1}{c_n^{p_0}} \int_E
\bigl|S_{\lceil
nt \rceil- \lceil ns \rceil} (f)\bigr|^{p_0} \,d\mu.
\end{eqnarray*}
Recall the notation $A_n^{(p_0)} =( \int_E |S_n(f) / a_n |^{p_0} \,d \mu_n
)^{1/p_0}$ in (\ref{eqexpressioncn}). We conclude that
\begin{eqnarray*}
&& P \Biggl( \Biggl|\sum_{k=1}^{\lceil nt \rceil}X_k^{(2)}
- \sum_{k=1}^{\lceil ns \rceil}X_k^{(2)}
\Biggr| \geq\lambda c_n \Biggr)
\\
&&\qquad \leq \frac{C_1+C_2}{ \lambda^{p_0}} \frac{1}{c_n^{p_0}}
\int_E \bigl|S_{\lceil
nt \rceil- \lceil ns \rceil} (f)\bigr|^{p_0} \,d\mu
\\
&&\qquad = \frac{C_1+C_2}{ \lambda^{p_0}}\frac{\mu(\varphi\leq\lceil nt
\rceil
- \lceil ns
\rceil)}{\mu(\varphi\leq n)} \biggl( \frac{a_{\lceil nt \rceil-
\lceil ns \rceil}}{a_n}
\biggr)^{p_0} \frac{(A_{\lceil
nt \rceil- \lceil ns \rceil}^{(
p_0)})^{p_0}}{c_n^{p_0} \mu(\varphi
\leq n)^{-1} a_n^{-p_0}}.
\end{eqnarray*}
It follows from (\ref{equnifintegral}) that
\[
\sup_{n \geq1, 0\leq s \leq t \leq L} A_{\lceil nt \rceil- \lceil ns
\rceil}^{( p_0)} < \infty.
\]
Next, we may, if necessary, increase $p_0$ in
(\ref{eqorginregularity}) to achieve $p_0>\alpha$. In that case,
the sequence $ c_n^{p_0} \mu(\varphi\leq n)^{-1}
a_n^{-p_0} \in RV_{(1-\beta)(p_0/\alpha-1)}$ diverges to
infinity, so for some constant $C_3>0$,
\[
\frac{1}{c_n^{p_0}} \int_E \bigl|S_{\lceil nt \rceil- \lceil ns \rceil}
(f)\bigr|^{p_0} \,d\mu\leq C_3 \frac{\mu(\varphi\leq\lceil n(t-s)
\rceil)}{\mu(\varphi\leq n)} \biggl(
\frac{a_{\lceil n(t-s)
\rceil}}{a_n} \biggr)^{p_0}.
\]
By the regular variation and the constraint $t-s\geq1/n$,
for every $0<\eta<\min(\beta,1-\beta)$, there is $C_4>0$, such that
\begin{eqnarray*}
\frac{\mu(\varphi\leq\lceil n(t-s) \rceil)}{\mu(\varphi\leq n)} &\leq& C_4 \biggl( \frac{\lceil n(t-s) \rceil}{n}
\biggr)^{1-\beta
-\eta} \leq2^{1-\beta-\eta} C_4 (t-s)
^{1-\beta-\eta},
\\
\frac{a_{\lceil n(t-s)\rceil}}{a_n} &\leq&2^{\beta-\eta} C_4 (t-s)
^{\beta-\eta}.
\end{eqnarray*}
Therefore, for some constant $C_5>0$,
\[
P \Biggl( \Biggl|\sum_{k=1}^{\lceil nt \rceil}X_k^{(2)}
- \sum_{k=1}^{\lceil ns \rceil}X_k^{(2)}
\Biggr| \geq\lambda c_n \Biggr) \leq C_5\frac{1}{\lambda^{p_0}}
(t-s)^{1+(p_0-1)\beta-(1+p_0)\eta}.
\]
Since $p_0> \alpha\geq1$, we can choose $\eta>0$ so small that
$1+(p_0-1)\beta- (1+p_0)\eta>0$. This establishes (\ref{epart2}).

Next, we take up the process $ ( X_n^{(1)} )$.
L\'evy--It{\^o} decomposition and the symmetry of the L\'evy measure
$\rho_1$ allow us to write, for any $K>0$,
\begin{eqnarray*}
\frac{1}{c_n} \sum_{k=1}^{\lceil nt \rceil}
X_k^{(1)} &\stackrel{d} {=}& \frac{1}{c_n} \sum
_{k=1}^{\lceil nt \rceil} \iint_{|x f_k|
\leq K
c_na_n^{-1}} x
f_k \,d\overline{N}_1 + \frac{1}{c_n} \sum
_{k=1}^{\lceil nt \rceil} \iint_{|x f_k| > Kc_n a_n^{-1}}  x
f_k \,dN_1
\\
&:=& Z_n^{(1,K)}(t) + Z_n^{(2,K)}(t),
\end{eqnarray*}
where $N_1$ and $\overline{N}_1$ are as above. Here, we first show that
or any $\epsilon>0$,
%
%
\begin{equation}
\label{elemma571} \lim_{K \to\infty} \limsup_{n \to\infty} P
\Bigl( \sup_{ 0 \leq t \leq L} \bigl|Z_n^{(2,K)}(t) \bigr| \geq
\epsilon \Bigr) =0.
\end{equation}
Consider first the case $1<\alpha<2$. Choose $0<\tau\leq2-\alpha$,
and define
\begin{eqnarray*}
\kappa(w) &=& \cases{ 1, &\quad if $0\leq w< 1$,
\cr
w^{-(\alpha+\tau)}, &\quad if
$w\geq1$,}
\\
g(w) &=& \bigl( (w+1)\kappa(w) \bigr)^{-1},\qquad w\geq0.
\end{eqnarray*}
Since $2g(w)/g(u)\geq1$ for $0\leq u\leq w$,
we have
\begin{eqnarray*}
&& P \Bigl( \sup_{ 0 \leq t \leq L} \bigl|Z_n^{(2,K)}(t) \bigr|
\geq\epsilon \Bigr)
\\
&&\qquad \leq P \Biggl( \iint_{\mathbb{R} \times E} |x| \sum
_{k=1}^{ nL } |f|\circ T^k \mathbf{1}
\bigl( |x||f|\circ T^k > Kc_n a_n^{-1}
\bigr) \,dN_1 \geq\epsilon c_n \Biggr)
\\
&&\qquad = P \Biggl( 2\iint_{\mathbb{R} \times E} |x| \sum_{k=1}^{ nL }
|f|\circ T^k g \bigl( |f|\circ T^k \bigr)
\frac{1}{g ( Kc_n
a_n^{-1}/|x| )} \,dN_1 \geq\epsilon c_n \Biggr)
\\
&&\qquad \leq \frac{2}{\epsilon} c_n^{-1} E \Biggl(
\iint_{\mathbb{R} \times E} |x| \sum_{k=1}^{ nL }
|f|\circ T^k g \bigl( |f|\circ T^k \bigr)
\frac{1}{g ( Kc_n
a_n^{-1}/|x| )} \,dN_1 \Biggr)
\\
&&\qquad \leq C_1 nc_n^{-1} \int_1^\infty
x \bigl( Kc_n a_n^{-1}/x+1 \bigr) \kappa \bigl(
Kc_n a_n^{-1}/x \bigr) \rho(dx),
\end{eqnarray*}
where $C_1>0$ is another constant. It is now straightforward to check
that for some constant $C_2>0$,
\[
\limsup_{n\to\infty} P \Bigl( \sup_{ 0 \leq t \leq L}
\bigl|Z_n^{(2,K)}(t) \bigr| \geq\epsilon \Bigr) \leq C_2
K^{-(\alpha-1)}.
\]
This implies (\ref{elemma571}).

On the other hand, let $0 < \alpha\leq1$. Recall that we are
assuming that the function $f$ is now bounded. We have
\begin{eqnarray*}
P \Bigl( \sup_{0 \leq t \leq L}\bigl|Z_n^{(2,K)}(t)\bigr| \geq
\epsilon \Bigr) &\leq& P \Bigl( \max_{k=1,\ldots,nL} N_1 \bigl
\{ (x,s)\dvtx  \bigl|xf_k(s)\bigr| > K c_n a_n^{-1}
\bigr\} \geq1 \Bigr)
\\
&\leq& E N_1 \Bigl\{ (x,s)\dvtx  |x| \max_{k=1,\ldots, nL}
|f_k| > K c_n a_n^{-1} \Bigr\}
\\
&=& 2 \int_E \rho_1 \biggl(
\frac{K c_n a_n^{-1}}{ \max_{k=1,\ldots, nL}
|f_k|}, \infty \biggr) \,d\mu.
\end{eqnarray*}
If we denote $\parallel
f \parallel= \sup_{x\in E} |f(x)|<\infty$, then we can use once again
Potter's bounds to see that for some constant $C_1>0$ and
$0<\xi<\alpha$,
\begin{eqnarray*}
&& \frac{\rho_1  (Kc_n a_n^{-1} (\max_k|f_k|)^{-1},\infty
)}{\rho_1(c_na_n^{-1},\infty)}
\\
&&\qquad  \leq C_1 \biggl( \biggl( \frac{1}{K} \max
_{k=1,\ldots,nL} |f_k| \biggr)^{\alpha-\xi} + \biggl(
\frac{1}{K} \max_{k=1,\ldots,nL} |f_k|
\biggr)^{\alpha+\xi} \biggr).
\end{eqnarray*}
Therefore by (\ref{eqwanderingmu}), (\ref{eqasycn}) and the fact
that $f$ is supported by $A$, for some constant $C_2 >0$,
\begin{eqnarray*}
&& P \Bigl( \sup_{0 \leq t \leq L}\bigl|Z_n^{(2,K)}(t)\bigr| \geq
\epsilon \Bigr)
\\
&&\qquad  \leq2C_1 \rho_1 \bigl(c_na_n^{-1},
\infty\bigr) \int_E \biggl( \frac{1}{K} \max
_{k=1,\ldots,nL} |f_k| \biggr)^{\alpha-\xi} + \biggl(
\frac{1}{K} \max_{k=1,\ldots,nL} |f_k|
\biggr)^{\alpha+\xi} \,d\mu
\\
&&\qquad \leq 2C_1 \rho_1 \bigl(c_n
a_n^{-1},\infty\bigr) \biggl( \biggl( \frac
{\parallel f
\parallel}{K}
\biggr)^{\alpha-\xi} + \biggl( \frac{\parallel f
\parallel
}{K} \biggr)^{\alpha+\xi}
\biggr) \mu(\varphi\leq nL)
\\
&&\qquad \leq C_2 \biggl( \biggl( \frac{\parallel f \parallel}{K} \biggr)^{\alpha-\xi}
+ \biggl( \frac{\parallel f \parallel}{K} \biggr)^{\alpha+\xi} \biggr)
\end{eqnarray*}
and (\ref{elemma571}) follows.

It remains to consider the processes $\{
Z_n^{(1,K)}(t), 0 \leq t \leq L \}$, $n=1,2,\ldots$ for a fixed
$K>0$. In the sequel, we drop the superscript $K$ for notational
convenience. We will show that exist $\gamma_1 > 1$,
and $B>0$ such that for all $0 \leq s < t \leq
L$, $n \geq1$, $t-s\geq1/n$ and $\lambda>0$,
%
%
\begin{equation}
P\bigl(\bigl|Z_n^{(1)}(t) - Z_n^{(1)}(s)\bigr|
\geq\lambda\bigr) \leq \frac{B}{\lambda^2} (t-s)^{\gamma_1}.\label{eqconcisegoal}
\end{equation}
Indeed, by Chebyshev's inequality and the fact that $f$ is supported
by $A$, we see that
\begin{eqnarray*}
&& P \bigl(\bigl|Z_n^{(1)}(t)-Z_n^{(1)}(s)\bigr|
\geq\lambda \bigr)
\\
&&\qquad \leq \frac{1}{\lambda^2 c_n^2} E \Biggl\llvert \sum
_{k=1}^{\lceil nt
\rceil-
\lceil ns \rceil} \iint_{|x f_k| \leq K c_na_n^{-1}} x
f_k \,d\overline{N}_1 \Biggr\rrvert ^2
\\
&&\qquad \leq\frac{2}{\lambda^2 c_n^2} \sum_{k=1}^{\lceil n(t-s) \rceil}
\sum_{l=1}^{\lceil n(t-s) \rceil} \int_E
|f_k f_l| \int
_0^{Kc_na_n^{-1} / |f_k| \vee|f_l|}
x^2 \rho_1(dx) \,d\mu.
\end{eqnarray*}
It follows from the Potter bounds and the fact that $\rho_1$ does not
assigns mass to the interval $(0,1)$ that for any $0<\xi<2-\alpha$
there is $C>0$ such that for all $a>0$ large enough and all $r>0$,
\[
\frac{\int_0^{ra}x^2  \rho_1(dx)}{\int_0^{a}x^2  \rho_1(dx)} \leq C \bigl( r^{2-\alpha-\xi} \vee r^{2-\alpha+\xi}
\bigr).
\]
Therefore, for all $n$ large enough, for some constant $C_1>0$,
\begin{eqnarray*}
&& P \bigl(\bigl|Z_n^{(1)}(t)-Z_n^{(1)}(s)\bigr|
\geq\lambda \bigr)
\\
&&\qquad \leq\frac{C_1}{\lambda^2 c_n^2} \sum_{k=1}^{\lceil n(t-s) \rceil}
\sum_{l=1}^{\lceil n(t-s) \rceil} \int_E
\frac{|f_k f_l|}{(|f_k| \vee|f_l|)^{2-\alpha-\xi}} \,d\mu \int_0^{c_na_n^{-1}}
x^2 \rho_1(dx)
\\
&&\quad\qquad{} +\frac{C_1}{\lambda^2 c_n^2} \sum_{k=1}^{\lceil n(t-s) \rceil}
\sum_{l=1}^{\lceil n(t-s) \rceil} \int_E
\frac{|f_k f_l|}{(|f_k| \vee|f_l|)^{2-\alpha+\xi}} \,d\mu \int_0^{c_na_n^{-1}}
x^2 \rho_1(dx).
\end{eqnarray*}
Note that by Karamata's theorem,
(\ref{eqwanderingmu}) and the definition (\ref{edefc}) of the
normalizing sequence $(c_n)$, there is $C_2>0$ such that
\[
\int_0^{c_na_n^{-1}} x^2
\rho_1(dx) \leq C_2 \frac{c_n^2}{na_n}.
\]
If $1<\alpha<2$, we impose also the constraint $\xi<\alpha-1$, and
use the relation
%
%
\begin{equation}
\label{eboundg1} \frac{|f_k f_l|}{(|f_k| \vee|f_l|)^{2-\alpha\pm\xi}} = \bigl( |f_k|\wedge|f_l|
\bigr) \bigl(|f_k| \vee |f_l| \bigr)^{\alpha-1\mp\xi},
\end{equation}
so that
\begin{eqnarray*}
&& \frac{1}{c_n^2} \sum_{k=1}^{\lceil n(t-s) \rceil} \sum
_{l=1}^{\lceil n(t-s) \rceil} \int_E
\frac{|f_k f_l|}{(|f_k| \vee|f_l|)^{2-\alpha\pm\xi}} \,d\mu \int_0^{c_na_n^{-1}}
x^2 \rho_1(dx)
\\
&&\qquad \leq C_2\frac{1}{na_n} \sum_{k=1}^{\lceil n(t-s) \rceil}
\sum_{l=1}^{\lceil n(t-s) \rceil} \int_E
\bigl( |f_k|\wedge|f_l| \bigr) \bigl(|f_k| \vee
|f_l| \bigr)^{\alpha-1\mp\xi} \,d\mu
\\
&&\qquad \leq2C_2 \frac{1}{na_n}
\\
&&\quad\qquad{}\times  \biggl[ \bigl\lceil n(t-s) \bigr\rceil\int
_E |f|^{\alpha\mp\xi} \,d\mu
\\
&&\hspace*{50pt}{}+ \sum_{k=1}^{\lceil n(t-s) \rceil-1} \sum
_{l=k+1}^{\lceil n(t-s) \rceil} \biggl( \int_E
|f_l| |f_k|^{\alpha
-1\mp\xi} \,d\mu + \int
_E |f_k||f_l|^{\alpha-1\mp\xi} \,d
\mu \biggr) \biggr]
\\
&&\qquad := J_n(1) + J_n(2)+J_n(3).
\end{eqnarray*}
The fact that $t-s>1/n$ and $(a_n)$ is regularly varying with the positive
exponent~$\beta$, shows that for any $1<\gamma_{1}<1+\beta$ there is
some constant $C_3>0$, such that for all $n=1,2,\ldots,$
\[
J_n(1)\leq C_3(t-s)^{\gamma_{1}}.
\]
Next, by the
duality relation (\ref{edualrel}),
\begin{eqnarray*}
J_n(2)&\leq&\frac{4C_2}{a_n}(t-s) \sum
_{k=1}^{\lceil n(t-s) \rceil} \int_E
|f_k||f|^{\alpha-1\mp\xi
} \,d\mu
\\
&=& \frac{4C_2}{a_n}(t-s) \int_A |f| \Biggl(\sum
_{k=1}^{\lceil n(t-s)
\rceil} \widehat{T}{}^k
|f|^{\alpha-1\mp\xi} \Biggr) \,d\mu.
\end{eqnarray*}
If $f$ is bounded, then by the
Darling--Kac property of the set $A$ we have, for some constants
$C_4,   C_5>0$,
\[
J_n(2)\leq C_4(t-s)\frac{a_{\lceil n(t-s) \rceil}}{a_n} \mu\bigl(|f|\bigr) \leq
C_5(t-s)^{\gamma_{1}},\qquad 1<\gamma_{1}<1+\beta
\]
by the regular variation of $(a_n)$. If, on the other hand, $A$ is a
uniform set for $|f|$, then we can write
\[
\sum_{k=1}^{\lceil n(t-s) \rceil} \widehat{T}{}^k
|f|^{\alpha-1\mp\xi} \leq \sum_{k=1}^{\lceil n(t-s) \rceil}
\widehat{T}{}^k \mathbf{1}_A + \sum
_{k=1}^{\lceil n(t-s) \rceil} \widehat{T}{}^k |f|
\]
and obtain the same bound on $J_2$ by using both the
Darling--Kac property and the uniform property
of the set $A$. A similar argument shows that, for some constant
$C_6>0$ we also have
\[
J_n(3)\leq C_6(t-s)^{\gamma_{1}},\qquad 1<
\gamma_{1}<1+\beta,
\]
which proves (\ref{eqconcisegoal}) in the case $1<\alpha<2$.

Finally, for $0< \alpha\leq1$ the same argument works, if we replace the
relation (\ref{eboundg1}) by
\begin{eqnarray*}
\frac{|f_k f_l|}{(|f_k| \vee|f_l|)^{1+\xi}}&\leq& \bigl( |f_k|\wedge |f_l|
\bigr)^{1-\xi},
\\
\frac{|f_k f_l|}{(|f_k| \vee
|f_l|)^{1-\xi}} &=& \bigl(|f_k|\wedge|f_l|
\bigr) \bigl(|f_k|\vee |f_l| \bigr)^\xi,
\end{eqnarray*}
respectively, if $\alpha=1$, and
\[
\frac{|f_k f_l|}{(|f_k| \vee|f_l|)^{2-\alpha\mp\xi}} \leq \bigl(|f_k|\wedge|f_l|
\bigr)^{\alpha\pm\xi}
\]
if $0<\alpha<1$.
This proves (\ref{eqconcisegoal}) in all cases, and hence, completes
the proof of the theorem.
\end{pf*}



%

\printaddresses


\begin{thebibliography}{50}

\bibitem[\protect\citeauthoryear{Aaronson}{1981}]{aaronson1981}
%
\begin{barticle}[mr]
\bauthor{\bsnm{Aaronson},~\bfnm{Jon}\binits{J.}}
(\byear{1981}).
\btitle{The asymptotic distributional behaviour of transformations
preserving infinite measures}.
\bjournal{J. Anal. Math.}
\bvolume{39}
\bpages{203--234}.
\bid{doi={10.1007/BF02803336}, issn={0021-7670}, mr={0632462}}
\end{barticle}
%
\bptok{imsref}%
\endbibitem

\bibitem[\protect\citeauthoryear{Aaronson}{1997}]{aaronson1997}
%
\begin{bbook}[mr]
\bauthor{\bsnm{Aaronson},~\bfnm{Jon}\binits{J.}}
(\byear{1997}).
\btitle{An Introduction to Infinite Ergodic Theory}.
\bseries{Mathematical Surveys and Monographs}
\bvolume{50}.
\bpublisher{Amer. Math. Soc.},
\blocation{Providence, RI}.
\bid{mr={1450400}}
\end{bbook}
%
\bptok{imsref}%
\endbibitem

\bibitem[\protect\citeauthoryear{Avram and Taqqu}{1992}]{avramtaqqu1992}
%
\begin{barticle}[mr]
\bauthor{\bsnm{Avram},~\bfnm{Florin}\binits{F.}} \AND
\bauthor{\bsnm{Taqqu},~\bfnm{Murad~S.}\binits{M.~S.}}
(\byear{1992}).
\btitle{Weak convergence of sums of moving averages in the {$\alpha
$}-stable domain of attraction}.
\bjournal{Ann. Probab.}
\bvolume{20}
\bpages{483--503}.
\bid{issn={0091-1798}, mr={1143432}}
\end{barticle}
%
\bptok{imsref}%
\endbibitem

\bibitem[\protect\citeauthoryear{Bertoin}{1996}]{bertoin1996}
%
\begin{bbook}[mr]
\bauthor{\bsnm{Bertoin},~\bfnm{Jean}\binits{J.}}
(\byear{1996}).
\btitle{L\'evy Processes}.
\bseries{Cambridge Tracts in Mathematics}
\bvolume{121}.
\bpublisher{Cambridge Univ. Press},
\blocation{Cambridge}.
\bid{mr={1406564}}
\end{bbook}
%
\bptok{imsref}%
\endbibitem

\bibitem[\protect\citeauthoryear{Billingsley}{1995}]{billingsley1995}
%
\begin{bbook}[mr]
\bauthor{\bsnm{Billingsley},~\bfnm{Patrick}\binits{P.}}
(\byear{1995}).
\btitle{Probability and Measure},
\bedition{3rd} ed.
\bpublisher{Wiley},
\blocation{New York}.
\bid{mr={1324786}}
\end{bbook}
%
\bptok{imsref}%
\endbibitem

\bibitem[\protect\citeauthoryear{Billingsley}{1999}]{billingsley1999}
%
\begin{bbook}[mr]
\bauthor{\bsnm{Billingsley},~\bfnm{Patrick}\binits{P.}}
(\byear{1999}).
\btitle{Convergence of Probability Measures},
\bedition{2nd} ed.
\bpublisher{Wiley},
\blocation{New York}.
\bid{doi={10.1002/9780470316962}, mr={1700749}}
\end{bbook}
%
\bptok{imsref}%
\endbibitem

\bibitem[\protect\citeauthoryear{Bingham}{1971}]{bingham1971}
%
\begin{barticle}[mr]
\bauthor{\bsnm{Bingham},~\bfnm{N.~H.}\binits{N.~H.}}
(\byear{1971}).
\btitle{Limit theorems for occupation times of {M}arkov processes}.
\bjournal{Z. Wahrsch. Verw. Gebiete}
\bvolume{17}
\bpages{1--22}.
\bid{mr={0281255}}
\end{barticle}
%
\bptok{imsref}%
\endbibitem

\bibitem[\protect\citeauthoryear{Cohen and
Samorodnitsky}{2006}]{cohensamorodnitsky2006}
%
\begin{barticle}[mr]
\bauthor{\bsnm{Cohen},~\bfnm{Serge}\binits{S.}} \AND
\bauthor{\bsnm{Samorodnitsky},~\bfnm{Gennady}\binits{G.}}
(\byear{2006}).
\btitle{Random rewards, fractional {B}rownian local times and stable
self-similar processes}.
\bjournal{Ann. Appl. Probab.}
\bvolume{16}
\bpages{1432--1461}.
\bid{doi={10.1214/105051606000000277}, issn={1050-5164}, mr={2260069}}
\end{barticle}
%
\bptok{imsref}%
\endbibitem

\bibitem[\protect\citeauthoryear{Darling and Kac}{1957}]{darlingkac1957}
%
\begin{barticle}[mr]
\bauthor{\bsnm{Darling},~\bfnm{D.~A.}\binits{D.~A.}} \AND
\bauthor{\bsnm{Kac},~\bfnm{M.}\binits{M.}}
(\byear{1957}).
\btitle{On occupation times for {M}arkoff processes}.
\bjournal{Trans. Amer. Math. Soc.}
\bvolume{84}
\bpages{444--458}.
\bid{issn={0002-9947}, mr={0084222}}
\end{barticle}
%
\bptok{imsref}%
\endbibitem

\bibitem[\protect\citeauthoryear{Davis and
Resnick}{1985}]{davisresnick1985}
%
\begin{barticle}[mr]
\bauthor{\bsnm{Davis},~\bfnm{Richard}\binits{R.}} \AND
\bauthor{\bsnm{Resnick},~\bfnm{Sidney}\binits{S.}}
(\byear{1985}).
\btitle{Limit theory for moving averages of random variables with
regularly varying tail probabilities}.
\bjournal{Ann. Probab.}
\bvolume{13}
\bpages{179--195}.
\bid{issn={0091-1798}, mr={0770636}}
\end{barticle}
%
\bptok{imsref}%
\endbibitem

\bibitem[\protect\citeauthoryear{Dobrushin and
Major}{1979}]{dobrushinmajor1979}
%
\begin{barticle}[mr]
\bauthor{\bsnm{Dobrushin},~\bfnm{R.~L.}\binits{R.~L.}} \AND
\bauthor{\bsnm{Major},~\bfnm{P.}\binits{P.}}
(\byear{1979}).
\btitle{Noncentral limit theorems for nonlinear functionals of
{G}aussian fields}.
\bjournal{Z. Wahrsch. Verw. Gebiete}
\bvolume{50}
\bpages{27--52}.
\bid{doi={10.1007/BF00535673}, issn={0044-3719}, mr={0550122}}
\end{barticle}
%
\bptok{imsref}%
\endbibitem

\bibitem[\protect\citeauthoryear{Dombry and
Guillotin-Plantard}{2009}]{dombryguillotin-plantard2009}
%
\begin{barticle}[mr]
\bauthor{\bsnm{Dombry},~\bfnm{C.}\binits{C.}} \AND
\bauthor{\bsnm{Guillotin-Plantard},~\bfnm{N.}\binits{N.}}
(\byear{2009}).
\btitle{Discrete approximation of a stable self-similar stationary
increments process}.
\bjournal{Bernoulli}
\bvolume{15}
\bpages{195--222}.
\bid{doi={10.3150/08-BEJ147}, issn={1350-7265}, mr={2546804}}
\end{barticle}
%
\bptok{imsref}%
\endbibitem

\bibitem[\protect\citeauthoryear{Ehm}{1981}]{ehm1981}
%
\begin{barticle}[mr]
\bauthor{\bsnm{Ehm},~\bfnm{W.}\binits{W.}}
(\byear{1981}).
\btitle{Sample function properties of multiparameter stable processes}.
\bjournal{Z. Wahrsch. Verw. Gebiete}
\bvolume{56}
\bpages{195--228}.
\bid{doi={10.1007/BF00535741}, issn={0044-3719}, mr={0618272}}
\end{barticle}
%
\bptok{imsref}%
\endbibitem

\bibitem[\protect\citeauthoryear{Getoor and
Kesten}{1972}]{getoorkesten1972}
%
\begin{barticle}[mr]
\bauthor{\bsnm{Getoor},~\bfnm{R.~K.}\binits{R.~K.}} \AND
\bauthor{\bsnm{Kesten},~\bfnm{H.}\binits{H.}}
(\byear{1972}).
\btitle{Continuity of local times for {M}arkov processes}.
\bjournal{Compos. Math.}
\bvolume{24}
\bpages{277--303}.
\bid{issn={0010-437X}, mr={0310977}}
\end{barticle}
%
\bptok{imsref}%
\endbibitem

\bibitem[\protect\citeauthoryear{Gradshteyn and
Ryzhik}{1994}]{gradshteynryzhik1994}
%
\begin{bbook}[mr]
\bauthor{\bsnm{Gradshteyn},~\bfnm{I.~S.}\binits{I.~S.}} \AND
\bauthor{\bsnm{Ryzhik},~\bfnm{I.~M.}\binits{I.~M.}}
(\byear{1994}).
\btitle{Table of Integrals, Series, and Products},
\bedition{5th} ed.
\bpublisher{Academic Press},
\blocation{Boston, MA}.
\bid{mr={1243179}}
\end{bbook}
%
\bptok{imsref}%
\endbibitem

\bibitem[\protect\citeauthoryear{Harris and
Robbins}{1953}]{harrisrobbins1953}
%
\begin{barticle}[mr]
\bauthor{\bsnm{Harris},~\bfnm{T.~E.}\binits{T.~E.}} \AND
\bauthor{\bsnm{Robbins},~\bfnm{Herbert}\binits{H.}}
(\byear{1953}).
\btitle{Ergodic theory of {M}arkov chains admitting an infinite
invariant measure}.
\bjournal{Proc. Natl. Acad. Sci. USA}
\bvolume{39}
\bpages{860--864}.
\bid{issn={0027-8424}, mr={0056873}}
\end{barticle}
%
\bptok{imsref}%
\endbibitem

\bibitem[\protect\citeauthoryear{Kallenberg}{2002}]{kallenberg2002}
%
\begin{bbook}[mr]
\bauthor{\bsnm{Kallenberg},~\bfnm{Olav}\binits{O.}}
(\byear{2002}).
\btitle{Foundations of Modern Probability},
\bedition{2nd} ed.
\bseries{Probability and Its Applications (New York)}.
\bpublisher{Springer},
\blocation{New York}.
\bid{mr={1876169}}
\end{bbook}
%
\bptok{imsref}%
\endbibitem

\bibitem[\protect\citeauthoryear{Kyprianou}{2006}]{kyprianou2006}
%
\begin{bbook}[mr]
\bauthor{\bsnm{Kyprianou},~\bfnm{Andreas~E.}\binits{A.~E.}}
(\byear{2006}).
\btitle{Introductory Lectures on Fluctuations of {L}\'evy Processes
with Applications}.
\bpublisher{Springer},
\blocation{Berlin}.
\bid{mr={2250061}}
\end{bbook}
%
\bptok{imsref}%
\endbibitem

\bibitem[\protect\citeauthoryear{Lamperti}{1962}]{lamperti1962}
%
\begin{barticle}[mr]
\bauthor{\bsnm{Lamperti},~\bfnm{John}\binits{J.}}
(\byear{1962}).
\btitle{Semi-stable stochastic processes}.
\bjournal{Trans. Amer. Math. Soc.}
\bvolume{104}
\bpages{62--78}.
\bid{issn={0002-9947}, mr={0138128}}
\end{barticle}
%
\bptok{imsref}%
\endbibitem

\bibitem[\protect\citeauthoryear{Maejima}{1983}]{maejima1983B}
%
\begin{barticle}[mr]
\bauthor{\bsnm{Maejima},~\bfnm{Makoto}\binits{M.}}
(\byear{1983}).
\btitle{On a class of self-similar processes}.
\bjournal{Z. Wahrsch. Verw. Gebiete}
\bvolume{62}
\bpages{235--245}.
\bid{doi={10.1007/BF00538799}, issn={0044-3719}, mr={0688988}}
\end{barticle}
%
\bptok{imsref}%
\endbibitem

\bibitem[\protect\citeauthoryear{Marcus and Rosen}{2006}]{marcusrosen2006}
%
\begin{bbook}[mr]
\bauthor{\bsnm{Marcus},~\bfnm{Michael~B.}\binits{M.~B.}} \AND
\bauthor{\bsnm{Rosen},~\bfnm{Jay}\binits{J.}}
(\byear{2006}).
\btitle{Markov Processes, {G}aussian Processes, and Local Times}.
\bseries{Cambridge Studies in Advanced Mathematics}
\bvolume{100}.
\bpublisher{Cambridge Univ. Press},
\blocation{Cambridge}.
\bid{doi={10.1017/CBO9780511617997}, mr={2250510}}
\end{bbook}
%
\bptok{imsref}%
\endbibitem

\bibitem[\protect\citeauthoryear{Maruyama}{1970}]{maruyama1970}
%
\begin{barticle}[auto]
\bauthor{\bsnm{Maruyama},~\bfnm{G.}\binits{G.}}
(\byear{1970}).
\btitle{Infinitely divisible processes}.
\bjournal{Theory Probab. Appl.}
\bvolume{15}
\bpages{3--22}.
\end{barticle}
%
\bptok{imsref}%
\endbibitem

\bibitem[\protect\citeauthoryear{Meerschaert and
Scheffler}{2004}]{meerschaertscheffler2004}
%
\begin{barticle}[mr]
\bauthor{\bsnm{Meerschaert},~\bfnm{Mark~M.}\binits{M.~M.}} \AND
\bauthor{\bsnm{Scheffler},~\bfnm{Hans-Peter}\binits{H.-P.}}
(\byear{2004}).
\btitle{Limit theorems for continuous-time random walks with infinite
mean waiting times}.
\bjournal{J. Appl. Probab.}
\bvolume{41}
\bpages{623--638}.
\bid{issn={0021-9002}, mr={2074812}}
\end{barticle}
%
\bptok{imsref}%
\endbibitem

\bibitem[\protect\citeauthoryear{Merlev{\`e}de, Peligrad and
Utev}{2006}]{merlevedepeligradutev2006}
%
\begin{barticle}[mr]
\bauthor{\bsnm{Merlev{\`e}de},~\bfnm{Florence}\binits{F.}},
\bauthor{\bsnm{Peligrad},~\bfnm{Magda}\binits{M.}} \AND
\bauthor{\bsnm{Utev},~\bfnm{Sergey}\binits{S.}}
(\byear{2006}).
\btitle{Recent advances in invariance principles for stationary sequences}.
\bjournal{Probab. Surv.}
\bvolume{3}
\bpages{1--36}.
\bid{doi={10.1214/154957806100000202}, issn={1549-5787}, mr={2206313}}
\end{barticle}
%
\bptok{imsref}%
\endbibitem

\bibitem[\protect\citeauthoryear{Paulauskas and
Surgailis}{2007}]{paulauskassurgailis2008}
%
\begin{barticle}[auto]
\bauthor{\bsnm{Paulauskas},~\bfnm{V.}\binits{V.}} \AND
\bauthor{\bsnm{Surgailis},~\bfnm{D.}\binits{D.}}
(\byear{2008}).
\btitle{On the rate of approximation in limit theorems for sums of
moving averages}.
\bjournal{Theory Probab. Appl.}
\bvolume{52}
\bpages{361--370}.
\end{barticle}
%
\bptok{imsref}%
\endbibitem

\bibitem[\protect\citeauthoryear{Pratt}{1960}]{pratt1960}
%
\begin{barticle}[mr]
\bauthor{\bsnm{Pratt},~\bfnm{John~W.}\binits{J.~W.}}
(\byear{1960}).
\btitle{On interchanging limits and integrals}.
\bjournal{Ann. Math. Statist.}
\bvolume{31}
\bpages{74--77}.
\bid{issn={0003-4851}, mr={0123673}}
\end{barticle}
%
\bptok{imsref}%
\endbibitem

\bibitem[\protect\citeauthoryear{Rajput and Rosi{\'n}ski}{1989}]{rajputrosinski1989}
%
\begin{barticle}[mr]
\bauthor{\bsnm{Rajput},~\bfnm{Balram~S.}\binits{B.~S.}} \AND
\bauthor{\bsnm{Rosi{\'n}ski},~\bfnm{Jan}\binits{J.}}
(\byear{1989}).
\btitle{Spectral representations of infinitely divisible processes}.
\bjournal{Probab. Theory Related Fields}
\bvolume{82}
\bpages{451--487}.
\bid{doi={10.1007/BF00339998}, issn={0178-8051}, mr={1001524}}
\end{barticle}
%
\bptok{imsref}%
\endbibitem

\bibitem[\protect\citeauthoryear{Resnick}{1987}]{resnick1987}
%
\begin{bbook}[mr]
\bauthor{\bsnm{Resnick},~\bfnm{Sidney~I.}\binits{S.~I.}}
(\byear{1987}).
\btitle{Extreme Values, Regular Variation, and Point Processes}.
\bpublisher{Springer},
\blocation{New York}.
\bid{mr={0900810}}
\end{bbook}
%
\bptok{imsref}%
\endbibitem

\bibitem[\protect\citeauthoryear{Resnick, Samorodnitsky and
Xue}{2000}]{resnicksamorodnitskyxue2000}
%
\begin{barticle}[mr]
\bauthor{\bsnm{Resnick},~\bfnm{Sidney}\binits{S.}},
\bauthor{\bsnm{Samorodnitsky},~\bfnm{Gennady}\binits{G.}} \AND
\bauthor{\bsnm{Xue},~\bfnm{Fang}\binits{F.}}
(\byear{2000}).
\btitle{Growth rates of sample covariances of stationary symmetric
{$\alpha$}-stable processes associated with null recurrent {M}arkov chains}.
\bjournal{Stochastic Process. Appl.}
\bvolume{85}
\bpages{321--339}.
\bid{doi={10.1016/S0304-4149(99)00081-2}, issn={0304-4149}, mr={1731029}}
\end{barticle}
%
\bptok{imsref}%
\endbibitem

\bibitem[\protect\citeauthoryear{Rosenblatt}{1956}]{rosenblatt1956}
%
\begin{barticle}[mr]
\bauthor{\bsnm{Rosenblatt},~\bfnm{M.}\binits{M.}}
(\byear{1956}).
\btitle{A central limit theorem and a strong mixing condition}.
\bjournal{Proc. Natl. Acad. Sci. USA}
\bvolume{42}
\bpages{43--47}.
\bid{issn={0027-8424}, mr={0074711}}
\end{barticle}
%
\bptok{imsref}%
\endbibitem

\bibitem[\protect\citeauthoryear{Rosi{\'n}ski}{1995}]{rosinski1995}
%
\begin{barticle}[mr]
\bauthor{\bsnm{Rosi{\'n}ski},~\bfnm{Jan}\binits{J.}}
(\byear{1995}).
\btitle{On the structure of stationary stable processes}.
\bjournal{Ann. Probab.}
\bvolume{23}
\bpages{1163--1187}.
\bid{issn={0091-1798}, mr={1349166}}
\end{barticle}
%
\bptok{imsref}%
\endbibitem

\bibitem[\protect\citeauthoryear{Rosi{\'n}ski and
Samorodnitsky}{1993}]{rosinskisamorodnitsky1993}
%
\begin{barticle}[mr]
\bauthor{\bsnm{Rosi{\'n}ski},~\bfnm{Jan}\binits{J.}} \AND
\bauthor{\bsnm{Samorodnitsky},~\bfnm{Gennady}\binits{G.}}
(\byear{1993}).
\btitle{Distributions of subadditive functionals of sample paths of
infinitely divisible processes}.
\bjournal{Ann. Probab.}
\bvolume{21}
\bpages{996--1014}.
\bid{issn={0091-1798}, mr={1217577}}
\end{barticle}
%
\bptok{imsref}%
\endbibitem

\bibitem[\protect\citeauthoryear{Rosi{\'n}ski and
Samorodnitsky}{1996}]{rosinskisamorodnitsky1996}
%
\begin{barticle}[mr]
\bauthor{\bsnm{Rosi{\'n}ski},~\bfnm{Jan}\binits{J.}} \AND
\bauthor{\bsnm{Samorodnitsky},~\bfnm{Gennady}\binits{G.}}
(\byear{1996}).
\btitle{Classes of mixing stable processes}.
\bjournal{Bernoulli}
\bvolume{2}
\bpages{365--377}.
\bid{doi={10.2307/3318419}, issn={1350-7265}, mr={1440274}}
\end{barticle}
%
\bptok{imsref}%
\endbibitem

\bibitem[\protect\citeauthoryear{Rosi{\'n}ski and {\.
Z}ak}{1996}]{rosinskizak1996}
%
\begin{barticle}[mr]
\bauthor{\bsnm{Rosi{\'n}ski},~\bfnm{Jan}\binits{J.}} \AND
\bauthor{\bsnm{{\.Z}ak},~\bfnm{Tomasz}\binits{T.}}
(\byear{1996}).
\btitle{Simple conditions for mixing of infinitely divisible processes}.
\bjournal{Stochastic Process. Appl.}
\bvolume{61}
\bpages{277--288}.
\bid{doi={10.1016/0304-4149(95)00083-6}, issn={0304-4149}, mr={1386177}}
\end{barticle}
%
\bptok{imsref}%
\endbibitem

\bibitem[\protect\citeauthoryear{Rosi{\'n}ski and {\.
Z}ak}{1997}]{rosinskizak1997}
%
\begin{barticle}[mr]
\bauthor{\bsnm{Rosi{\'n}ski},~\bfnm{Jan}\binits{J.}} \AND
\bauthor{\bsnm{{\.Z}ak},~\bfnm{Tomasz}\binits{T.}}
(\byear{1997}).
\btitle{The equivalence of ergodicity of weak mixing for infinitely
divisible processes}.
\bjournal{J. Theoret. Probab.}
\bvolume{10}
\bpages{73--86}.
\bid{doi={10.1023/A:1022690230759}, issn={0894-9840}, mr={1432616}}
\end{barticle}
%
\bptok{imsref}%
\endbibitem

\bibitem[\protect\citeauthoryear{Roy}{2007}]{roy2008}
%
\begin{barticle}[mr]
\bauthor{\bsnm{Roy},~\bfnm{Emmanuel}\binits{E.}}
(\byear{2007}).
\btitle{Ergodic properties of {P}oissonian {ID} processes}.
\bjournal{Ann. Probab.}
\bvolume{35}
\bpages{551--576}.
\bid{doi={10.1214/009117906000000692}, issn={0091-1798}, mr={2308588}}
\end{barticle}
%
\bptok{imsref}%
\endbibitem

\bibitem[\protect\citeauthoryear{Samorodnitsky}{2004}]{samorodnitsky2004a}
%
\begin{barticle}[mr]
\bauthor{\bsnm{Samorodnitsky},~\bfnm{Gennady}\binits{G.}}
(\byear{2004}).
\btitle{Extreme value theory, ergodic theory and the boundary between
short memory and long memory for stationary stable processes}.
\bjournal{Ann. Probab.}
\bvolume{32}
\bpages{1438--1468}.
\bid{doi={10.1214/009117904000000261}, issn={0091-1798}, mr={2060304}}
\end{barticle}
%
\bptok{imsref}%
\endbibitem

\bibitem[\protect\citeauthoryear{Samorodnitsky}{2005}]{samorodnitsky2005a}
%
\begin{barticle}[mr]
\bauthor{\bsnm{Samorodnitsky},~\bfnm{Gennady}\binits{G.}}
(\byear{2005}).
\btitle{Null flows, positive flows and the structure of stationary
symmetric stable processes}.
\bjournal{Ann. Probab.}
\bvolume{33}
\bpages{1782--1803}.
\bid{doi={10.1214/009117905000000305}, issn={0091-1798}, mr={2165579}}
\end{barticle}
%
\bptok{imsref}%
\endbibitem

\bibitem[\protect\citeauthoryear{Samorodnitsky}{2006}]{samorodnitsky2006LRD}
%
\begin{barticle}[mr]
\bauthor{\bsnm{Samorodnitsky},~\bfnm{Gennady}\binits{G.}}
(\byear{2006}).
\btitle{Long range dependence}.
\bjournal{Found. Trends Stoch. Syst.}
\bvolume{1}
\bpages{163--257}.
\bid{doi={10.1561/0900000004}, issn={1551-3106}, mr={2379935}}
\end{barticle}
%
\bptok{imsref}%
\endbibitem

\bibitem[\protect\citeauthoryear{Samorodnitsky and
Taqqu}{1994}]{samorodnitskytaqqu1994}
%
\begin{bbook}[mr]
\bauthor{\bsnm{Samorodnitsky},~\bfnm{Gennady}\binits{G.}} \AND
\bauthor{\bsnm{Taqqu},~\bfnm{Murad~S.}\binits{M.~S.}}
(\byear{1994}).
\btitle{Stable Non-{G}aussian Random Processes}.
\bpublisher{Chapman \& Hall},
\blocation{New York}.
\bid{mr={1280932}}
\end{bbook}
%
\bptok{imsref}%
\endbibitem

\bibitem[\protect\citeauthoryear{Sato}{1999}]{sato1999}
%
\begin{bbook}[mr]
\bauthor{\bsnm{Sato},~\bfnm{Ken-iti}\binits{K.-i.}}
(\byear{1999}).
\btitle{L\'evy Processes and Infinitely Divisible Distributions}.
\bpublisher{Cambridge Univ. Press},
\blocation{Cambridge}.
\bid{mr={1739520}}
\end{bbook}
%
\bptok{imsref}%
\endbibitem

\bibitem[\protect\citeauthoryear{Taqqu}{1979}]{taqqu1979}
%
\begin{barticle}[mr]
\bauthor{\bsnm{Taqqu},~\bfnm{Murad~S.}\binits{M.~S.}}
(\byear{1979}).
\btitle{Convergence of integrated processes of arbitrary {H}ermite rank}.
\bjournal{Z. Wahrsch. Verw. Gebiete}
\bvolume{50}
\bpages{53--83}.
\bid{doi={10.1007/BF00535674}, issn={0044-3719}, mr={0550123}}
\end{barticle}
%
\bptok{imsref}%
\endbibitem

\bibitem[\protect\citeauthoryear{Thaler}{2001}]{thaler2001}
%
\begin{bmisc}[auto:STB|2014/01/06|10:16:28]
\bauthor{\bsnm{Thaler},~\bfnm{M.}\binits{M.}}
(\byear{2001}).
\bhowpublished{Infinite ergodic theory. The Dynamic Odyssey course, CIRM.}
\end{bmisc}
%
\bptok{imsref}%
\endbibitem

\bibitem[\protect\citeauthoryear{Thaler and Zweim{\"
u}ller}{2006}]{thalerzweimuller2006}
%
\begin{barticle}[mr]
\bauthor{\bsnm{Thaler},~\bfnm{M.}\binits{M.}} \AND
\bauthor{\bsnm{Zweim{\"u}ller},~\bfnm{R.}\binits{R.}}
(\byear{2006}).
\btitle{Distributional limit theorems in infinite ergodic theory}.
\bjournal{Probab. Theory Related Fields}
\bvolume{135}
\bpages{15--52}.
\bid{doi={10.1007/s00440-005-0454-3}, issn={0178-8051}, mr={2214150}}
\end{barticle}
%
\bptok{imsref}%
\endbibitem

\bibitem[\protect\citeauthoryear{Whitt}{2002}]{whitt2002}
%
\begin{bbook}[mr]
\bauthor{\bsnm{Whitt},~\bfnm{Ward}\binits{W.}}
(\byear{2002}).
\btitle{Stochastic-Process Limits: An Introduction to
Stochastic-Process Limits and Their Application to Queues}.
\bpublisher{Springer},
\blocation{New York}.
\bid{mr={1876437}}
\end{bbook}
%
\bptok{imsref}%
\endbibitem

\bibitem[\protect\citeauthoryear{Zolotarev}{1986}]{zolotarev1986}
%
\begin{bbook}[mr]
\bauthor{\bsnm{Zolotarev},~\bfnm{V.~M.}\binits{V.~M.}}
(\byear{1986}).
\btitle{One-Dimensional Stable Distributions}.
\bpublisher{Amer. Math. Soc.},
\blocation{Providence, RI}.
\bnote{Translated from the Russian by H. H. McFaden, Translation
edited by Ben Silver}.
\bid{mr={0854867}}
\end{bbook}
%
\bptok{imsref}%
\endbibitem

\bibitem[\protect\citeauthoryear{Zweim{\"u}ller}{2000}]{zweimuller2000}
%
\begin{barticle}[mr]
\bauthor{\bsnm{Zweim{\"u}ller},~\bfnm{Roland}\binits{R.}}
(\byear{2000}).
\btitle{Ergodic properties of infinite measure-preserving interval
maps with indifferent fixed points}.
\bjournal{Ergodic Theory Dynam. Systems}
\bvolume{20}
\bpages{1519--1549}.
\bid{doi={10.1017/S0143385700000821}, issn={0143-3857}, mr={1786727}}
\end{barticle}
%
\bptok{imsref}%
\endbibitem

\bibitem[\protect\citeauthoryear{Zweim{\"u}ller}{2007a}]{zweimuller2007}
%
\begin{barticle}[mr]
\bauthor{\bsnm{Zweim{\"u}ller},~\bfnm{Roland}\binits{R.}}
(\byear{2007}a).
\btitle{Infinite measure preserving transformations with compact first
regeneration}.
\bjournal{J. Anal. Math.}
\bvolume{103}
\bpages{93--131}.
\bid{doi={10.1007/s11854-008-0003-y}, issn={0021-7670}, mr={2373265}}
\end{barticle}
%
\bptok{imsref}%
\endbibitem

\bibitem[\protect\citeauthoryear{Zweim{\"u}ller}{2007b}]{zweimuller2007a}
%
\begin{barticle}[mr]
\bauthor{\bsnm{Zweim{\"u}ller},~\bfnm{Roland}\binits{R.}}
(\byear{2007}b).
\btitle{Mixing limit theorems for ergodic transformations}.
\bjournal{J. Theoret. Probab.}
\bvolume{20}
\bpages{1059--1071}.
\bid{doi={10.1007/s10959-007-0085-y}, issn={0894-9840}, mr={2359068}}
\end{barticle}
%
\bptok{imsref}%
\endbibitem

\bibitem[\protect\citeauthoryear{Zweim{\"u}ller}{2009}]{zweimuller2009}
%
\begin{bmisc}[auto:STB|2014/01/06|10:16:28]
\bauthor{\bsnm{Zweim{\"u}ller},~\bfnm{R.}\binits{R.}}
(\byear{2009}).
\bhowpublished{Surrey notes on infinite ergodic theory. Lecture notes,
Surrey Univ.}
\end{bmisc}
%
\bptok{imsref}%
\endbibitem

\end{thebibliography}
\end{document}